\documentclass[preprint]{imsart}

\usepackage[T1]{fontenc}
\usepackage[utf8]{inputenc}
\usepackage[english]{babel}
\usepackage{amsmath,amsthm,amssymb,amsfonts,epic,latexsym,hyperref,bbm}
\usepackage{graphicx}
\usepackage{cancel}
\usepackage{rotating}
\usepackage{dsfont}
\usepackage{geometry}
\usepackage{stmaryrd}
\usepackage{ifthen}
\usepackage{verbatim}
\usepackage{tabularx}
\usepackage{framed}
\usepackage{longtable}
\usepackage{mathtools}
\usepackage[shortlabels]{enumitem}
\usepackage[usenames, dvipsnames]{xcolor}
\usepackage{eurosym}
\usepackage{tikz}

\numberwithin{equation}{section}

\definecolor{mycolorred}{rgb}{1, 0, 0}

\definecolor{mycolorblue}{rgb}{0, 0, 1}

\definecolor{mycolorpink}{rgb}{1, 0, 1}

\newtheorem{theorem}{Theorem}[section]
\newtheorem{proposition}[theorem]{Proposition}
\newtheorem{corollary}[theorem]{Corollary}
\newtheorem{lemma}[theorem]{Lemma}
\newtheorem{remark}[theorem]{Remark}

\newtheorem{assumption}[theorem]{Assumption}
\newtheorem{example}[theorem]{Example}

\def\<{\langle}
\def\>{\rangle}

\def\E{{\mathbb E}}
\def\R{{\mathbb R}}

\def\N{{\mathbb N}}

\def\MR #1{\textrm{MR }{#1}}

\begin{document}

\begin{frontmatter}

\title{Conditional propagation of chaos for systems of particles with polynomial interactions}

\runtitle{Particle systems with polynomial interactions}

\begin{aug}

\author{\fnms{Vlad} \snm{Bally}\thanksref{m1}\ead[label=e1]{vlad.bally@univ-eiffel.fr}}
\author{\fnms{Eva} \snm{L\"ocherbach}\thanksref{m2}\ead[label=e2]{eva.loecherbach@polytechnique.edu}}

\address{\thanksmark{m1}LAMA, Universit\'e Gustave Eiffel, and MATHRISK, INRIA,  France
  \thanksmark{m2}CMAP, Ecole Polytechnique, Institut Polytechnique de Paris, France
  }

\runauthor{V. Bally and E. L\"ocherbach}

 \end{aug}

\begin{abstract}
We consider $N$ interacting particles taking values in $ \R $ and evolving in discrete time with time step $ 1/n.$ Each particle follows an individual
dynamics, driven by its own noise, which depends on the total configuration via a polynomial interaction. Moreover, all particles receive the same collective
perturbation. This perturbation is given by a multilinear homogeneous sum of order $ m \geq 2$ of independent centered random fields which are evaluated at the 
current positions of the particles. We show that, as $N \to \infty $ and $ n \to \infty, $ the system converges to the solution of a conditional McKean-Vlasov equation
which is driven by an additional Brownian motion which emerges in the limit from the collective perturbations. Moreover we give explicit bounds for the 
Wasserstein-$2-$distance. Our proof combines the invariance principle of Mossel, O'Donnell and Oleszkiewicz \cite{MDO}, the CLT for Wiener chaoses 
in the Wasserstein-$2-$form of Ledoux, Nourdin 
and Peccati \cite{LNP} and a conditional law of large numbers for weighted $U-$statistics proven in this paper.

\end{abstract}


\end{frontmatter}

\section{Introduction}

\label{sec:1} In this paper we establish a dynamical invariance principle for
 a system of $N$ interacting real valued particles and show how this gives rise to the creation of some common noise
in the limit. In our model,  particles evolve in discrete time, and they interact through
polynomial interactions which are described below.  We will show that, as $N \to \infty$ and as the time
step converges to $0,$ 
the system can be described by a limit stochastic differential equation
which is of conditional McKean-Vlasov type, driven by some limit common
noise. This limit common noise is a Brownian motion that arises as the
consequence of a dynamical version of the invariance principle of Mossel et al. \cite{MDO}
and the CLT for multiple stochastic integrals of Nualart and Peccati \cite{NuPe}, Nourdin and Peccati \cite{NP} and of Ledoux, Nourdin and Peccati \cite{LNP}.

We now describe the particle system more precisely. For fixed $N,$ the
number of particles, and $m,$ the order of interactions, we denote by $%
\Lambda _{N,m}$ the multi-indices $\alpha =(\alpha _{1},...,\alpha _{m})\in
\{1,...,N \}^{m}.$ We consider a collection of non-negative  coefficients $%
c_{N}(\alpha ) ,\alpha \in \Lambda _{N,m}, $ which are
symmetric (that is, invariant under permutation of its arguments) and null
on the diagonals. 

Moreover we consider an i.i.d. sequence of random fields $v_{k}^{i}:\Omega \times {%
\mathbb{R}}\rightarrow {\mathbb{R}}, k\in {\mathbb{N}},i=1,...,N,$ which verify, for each fixed $x\in {\mathbb{R}} , $%
\begin{equation}  \label{eq:vik}
E(v_{k}^{i}(x))=0,\quad E(v_{k}^{i}(x)^{2})=f^{2}(x) ,
\end{equation}%
where $f$ is a positive, bounded Lipschitz continuous function which is bounded away from $0.$  We also assume that for some $q\geq 3$ we have that 
\begin{equation}  \label{eq:Mq}
\sup_{k,i}\sup_{x}E\left\vert v_{k}^{i}(x)\right\vert ^{q}<\infty .
\end{equation}

We associate to these random fields the multi-linear homogeneous sum 
\begin{equation}  \label{eq:model}
S_{k}(x)=\sum_{\alpha \in \Lambda _{N,m}}c_{N}(\alpha
)\prod_{i=1}^{m}v_{k}^{\alpha _{i}}(x^{\alpha _{i}}), \; x=(x^{1},...,x^{N})
\in {\mathbb{R}}^N .
\end{equation}
We interpret this quantity $S_k( x) $ as a bombardment term to which each particle is equally exposed during one time step of its
evolution. So in our system the random input received by a particle is not a sum of independent individual contributions, but is produced by groups of $m$ particles acting together, as in models with higher-order (or multi-body) interactions. The homogeneous sum $S_k$ is the simplest model of such a collective input: each group $\alpha$ of $m$ particles contributes the product of their individual random signals, with weight $c_N(\alpha)$, and for $m = 1$ we recover the linear bombardment term considered in 
\cite{EvaVlad}. In this paper we consider $ m \geq 2.$ 

Moreover we consider another order of interactions $ \bar m \geq 1 $ and another collection of non-negative  coefficients $%
\bar c_{N}(\alpha )\in {\mathbb{R}},\alpha \in \Lambda _{N, \bar m},$ together with some Lipschitz continuous coefficients $b$ and $a$ and a 
bounded Lipschitz continuous
function $\varphi .$ Then we define 
\begin{equation*}
\mathbb{T}(x)=\sum_{\alpha \in \Lambda _{N,\bar m}}\bar c_{N}^{2}(\alpha
)\prod_{i=1}^{\bar m}\varphi (x^{\alpha _{i}}) .
\end{equation*}

These multi-linear homogeneous sums are the main building blocks of the
particle system which we introduce now. 

We fix some initial configuration $%
X_{0}=(X_{0}^{1},...,X_{0}^{N})$ such that the $ X_0^i , 1 \le i \le N, $ are i.i.d. and square integrable. 
We also consider a sequence $\Delta
_{k}^{i},i=1,...,N,$ $k\in \N ,$ of independent standard normal random
variables. We suppose that the families $ (v^i_k), (\Delta_k^i ) $ and $ (X_0^i ) $ are independent. 

Putting  $t_{k}=\frac{k}{n},$ $ k \geq 0,$ 
we then construct the chain 
\begin{equation}
X_{t_{k+1}}^{n,i}=X_{t_k}^{n,i}+b(X_{t_k}^{n,i},\mathbb{T}(X_{t_k}^{n}))%
\frac{1}{n}+a(X_{t_k}^{n,i},\mathbb{T}(X_{t_k}^{n}))\frac{1}{\sqrt{n}}%
\Delta _{k}^{i}+\frac{1}{\sqrt{n}}S_{k}(X_{t_k}^{n}), k \geq 0, 
\label{I1}
\end{equation}%
$ 1 \le i \le N, $ 
for some fixed time step $1/n,$ with $n\in {\mathbb{N}}.$ 

We stress that the above system is in general not exchangeable and that the interactions are not of mean field type and do not depend on the empirical measure of the system, since 
they are of polynomial order. Finally, 
notice that every particle $i$ receives the same random and centered input 
$\frac{1}{\sqrt{n}}S_{k}(X_{t_k}^{n}) .$ This will give rise to the creation of a common noise in the limit system.

We will show that, as both $N\rightarrow \infty $ and $1/n\rightarrow 0,$
this evolution converges to the limit stochastic differential equation
which is given by 
\begin{multline}
\mathbb{X}_{t}^{i} =X_{0}^{i}+\int_{0}^{t}b(\mathbb{X}_{s}^{i},(\E_{B}%
\varphi (\mathbb{X}_{s}^{i}))^{\bar m})ds+\int_{0}^{t}a(\mathbb{X}%
_{s}^{i},(\E_{B}\varphi (\mathbb{X}_{s}^{i}))^{\bar m})dW_{s}^{i}  \label{eq:SDE}
\\
+\int_{0}^{t}(\E_{B}f^{2}(\mathbb{X}_{s}^{i}))^{m/2}dB_{s},\quad i \geq 1.
\end{multline}%
In the above equation, $B$ is a standard Brownian motion, representing a source of common noise for 
all particles, and $\E_{B}$ denotes the
conditional expectation with respect to $\sigma \{B_{s},s\geq 0\}.$ Moreover, $%
W^{i},i \geq 1, $ are individual independent Brownian motions, which are
independent of $B.$ Finally, the initial conditions $ X_0^i, i \geq 1, $ are i.i.d. square integrable 
random variables which 
are independent of $ W^j, j \geq 1, $ and of $ B.$ 

We will work under standard Lipschitz and boundedness assumptions implying that
equation (\ref{eq:SDE}) possesses a unique strong solution. Notice that its structure is the one of 
a conditional McKean–Vlasov equation. Since the initial values are i.i.d., for every $t>0,$ the $\mathbb X^i_t$, $i\ge1$, are i.i.d. 
conditionally on $B$; this is the conditional propagation of chaos property which we will establish for the finite system.

A main motivation of this work is to understand what becomes of the invariance principle of Mossel et al. \cite{MDO}
when it is superposed on a stochastic dynamics, that is, when the multilinear sums $S_k$ are evaluated, at each time step, along the random configuration $X^n_{t_k}$ which they themselves help to produce. Looking at it this way, our main result can thus be read as a dynamical invariance principle. 

\subsection*{A strong error bound for the conditional propagation of chaos}
We follow a strong approach to quantify the convergence of the finite particle system to the limit system. 
We will produce an explicit coupling of the finite system with a discrete time version of the limit system, following a general strategy that
has been given in \cite{EvaVlad}. Doing so, our main result is a strong error bound between the finite particle system
and its associated limit, for fixed $ N $ and $n.$ This bound  is obtained under suitable assumptions on
the collection $c_N (\alpha)$ and $ \bar c_N (\alpha),$ which are stated precisely in Section \ref{sec:2} below. We work with the Wasserstein distance $\mathcal{W}_{2}$ of
probability measures on ${\mathbb{R}}^{N},$ where ${\mathbb{R}}^N$ is
equipped with the norm 
\begin{equation}  \label{eq:distanceRN}
\left\vert x\right\vert _{N}^{2}=\frac{1}{N}\sum_{i=1}^{N}\left\vert
x^{i}\right\vert ^{2}.
\end{equation}

\begin{theorem}[Preliminary version]
\label{theo:main0} Suppose that (\ref{eq:Mq}) holds for every $q\in {\mathbb{N}}$ and that the initial values $X_{0}^{1},...,X_{0}^{N}$ are i.i.d.
and square integrable. 
Then, under suitable assumptions on $c_{N}$ and on $ \bar c_N$ (see (\ref{McKV6}) below), for every $p<\frac{1}{3},$ 
\begin{equation}
\mathcal{W}_{2}^{2}(X_{1}^{n},\mathbb{X}_{1})\leq C_{p}(\frac{1}{N^{p/2 }}%
+ \sum_{r=1}^{ m-1} \| c_N \otimes_r c_N\|_{m-r, 2 }^2 +\frac{1}{n}),  \label{M1}
\end{equation}%
where the contraction factors $ \| c_N \otimes_r c_N\|_{m-r, 2 }$ are defined in (\ref{inf4}) below. In the above equation, 
$\mathbb{X}_{1}$ denotes the first $N$ components of $ \mathbb{X} $ evaluated at time $1.$ 
\end{theorem}

We stress that this result holds true for any possible choice of $N$ and $n.$ 

In our paper, we work within a more general framework than the one considered above. The precise statement of our main result within this general framework is given in Theorem \ref{theo:main} below, and the proof of this result is given in Sections \ref{sec:3} and \ref%
{sec:4} below.

\begin{remark}
The restriction $ p < 1/3$ is due to the fact that at some point in our proof we have to replace a Zolotarev type distance $d_3$ by a Wasserstein-$2-$distance, see {\bf Step 2} and also Lemma \ref{lem:21} below.
\end{remark}

\subsection*{The proof strategy}
Recall that  $t_{k}=\frac{k}{n}.$ We will show that we can interpret $X^n_{t_k} $ as
approximation of an Euler scheme, evaluated at time $t_k, $ of the limit SDE 
$\mathbb{X} .$ In what follows we wish to give the main ideas of our proof and we decided to concentrate on the creation of the common noise. 
Therefore, in what follows, all estimates will be given for the particular case where $ a = b = 0.$ 
We refer to our main result, Theorem \ref{theo:main}, for the general statement. 

\noindent \textbf{Step 1 }We fix $k$ and study what happens in the time
interval $[t_{k},t_{k+1}).$ This is done in several steps.

\noindent \textbf{Step 1.a (Invariance principle) } Let $%
G_{l}^{i},i=1,...,N,l\geq 1,$ be independent standard normal random
variables and define 
\begin{equation*}
S_{k}(x,G) =\sum_{\alpha \in \Lambda _{N,m}}c_{N}(\alpha
)\prod_{i=1}^{m}f(x^{\alpha _{i}})\prod_{i=1}^{m}G_{k}^{\alpha _{i}}.
\end{equation*}
Consider, for any two real valued random variables $U,V,$ the Zolotarev type
distance%
\begin{equation}  \label{eq:d3}
d_{3}(U,V)=\sup \{\left\vert E(h (U))-E(h (V))\right\vert
:\sum_{l=0}^{3}\left\Vert h ^{(l)}\right\Vert _{\infty }\leq 1\} .
\end{equation}
A first main ingredient of our approach is to use the invariance principle of Mossel et al. \cite{MDO}
which we apply in a framework where the multi-linear homogeneous sums depend on the configuration $x$.
More precisely, we introduce $ Z_k^i (x^i) = \frac{v^i_k ( x^i) }{f(x^i ) }, $ and we write $ S_k ( x) = \sum_\alpha c_N( \alpha) \prod_{i=1}^m f ( x^{\alpha_i}) \prod_{i=1}^m Z_k^{ \alpha_i} (x^{\alpha_i} ) .$ 
Then the invariance principle allows to replace $ Z_k^i (x^i) $ by $G_k^i, $ and so we obtain
\begin{equation}
d_{3}(S_{k}(x),S_{k}(x,G)) \leq C \delta _{N}(c_{N}) .  \label{I2}
\end{equation}%
Here, the \textit{influence factor} $\delta _{N}(c_{N})$ is defined by 
\begin{equation}
\delta _{N}^{2}(c_{N})=\max_{i=1,...,N}\sum_{\alpha \in \Lambda
_{N,m}^{i}}c_{N}^{2}(\alpha ) ,  \label{inf1}
\end{equation}%
and $\Lambda _{N,m}^{i}$ represents the set of multi-indices $\alpha \in
\Lambda _{N,m}$ which contain $i.$

We will always work under the hypothesis that $\delta_N (c_N) , \delta_N ( \bar c_N) \le  C N^{ - 1/2}  $ 
for some fixed constant $ C > 0.$ 

\noindent \textbf{Step 1.b (CLT) }We notice that $S_{k}(x,G)$ is actually a
multiple stochastic integral of order $m$ with piecewise constant
coefficients. As a consequence, we may use the CLT for multiple stochastic
integrals due to Nourdin and Peccati \cite{NP} and to Ledoux, Nourdin and Peccati \cite{LNP}. It implies that 
\begin{equation}
{W}_{2}^{2}(S_{k}(x,G),\sigma_N (x)\Delta )\leq C {\varepsilon 
}_{N}(c_N),  \label{I3}
\end{equation}%
where, for a constant $C_m$ depending on $m,$ 
$$ {\varepsilon }_{N}(c_N) = C_m \sum_{r=1}^{ m-1} \| c_N \otimes_r c_N\|_{m-r, 2 }^2, $$
is a quantity based on the
contraction operators on the Wiener space which are introduced in (\ref{inf4}) below. In the above formula, $\Delta $ is a standard normal
random variable, $W_2$ is the Wasserstein $2-$distance on $ \R,$ and 
\begin{equation}\label{eq:sigma}
\sigma_N ^2 (x ) := m! \sum_{\alpha
}c_{N}^{2}(\alpha )\prod_{j=1}^{m}f^2(x^{\alpha _{j}}).
\end{equation}

\begin{remark}
So far, in \textbf{Step 1.a } and \textbf{Step 1.b }, we have used the
same strategy as in \cite{NPR}. The only difference is that our coefficients
depend on the configuration $x.$ This is specific for our problem and will
be important in the sequel.
\end{remark}

\noindent \textbf{Step 2 (Concatenation)} Recall that  $t_{k}=\frac{k}{n}.$ We define now the chain 
\begin{multline}\label{eq:introdbarX}
\overline{X}_{t_{k+1}}^{n,i}=\overline{X}_{t_{k}}^{n,i}+\frac1n b(\overline X_{t_k}^{n,i},%
\mathbb{T}(\overline X_{t_k}^{n})) +\frac{1}{\sqrt{n}} a(\overline X_{t_k}^{n,i},\mathbb{T}(\overline X_{t_k}^{n}%
))\Delta _{k}^{i}\\
+\sigma_N (\overline{X}%
_{t_{k}}^{n })(B_{t_{k+1}}-B_{t_{k}}),0\leq k<n, 1 \le i \le N, 
\end{multline}%
where $B$ is a standard one-dimensional Brownian motion. We would like to
iterate the estimates (\ref{I2}) and (\ref{I3}). Since $%
B_{t_{k+1}}-B_{t_{k}}\sim (t_{k+1}-t_{k})^{1/2}\Delta ,$ at each step we
have a dependency on the time increment given by $(t_{k+1}-t_{k})^{1/2}$
which in our case equals $n^{-1/2}.$ We have in total $n$ such time
intervals to recover $[0,1],$ so this contribution clearly explodes as $%
n\rightarrow \infty .$

To tackle this problem, we replace the distance $d_{3}$\ (which appears in (%
\ref{I2})) by the Wasserstein distance $W_{2}.$ To do so, we rely on Lemma
B2 in \cite{ABCK}. Once this change is done, we may use a coupling argument
in order to transform the computations in $W_{2}$ into computations in $L^{2}
$ - and then we use standard techniques for discrete time martingales in order
to replace $(t_{k+1}-t_{k})^{1/2}$ by $(t_{k+1}-t_{k}).$ The precise
arguments are gathered in Proposition 4.2 of \cite{EvaVlad}. We stress that
this step is rather delicate and has to be done carefully.

As a consequence of these arguments, we then arrive at a control 
\begin{equation*}
\mathcal{W}_{2}^{2}(X_{t_{k}}^{n},\overline{X}_{t_{k}}^{n})\leq C(t_{k} - t_0) (
\delta^p_{N}(c) +  \varepsilon_N ( c_N) ).
\end{equation*}

\noindent \textbf{Step 3} We now want to replace the variance $\sigma_N^{2}(x)
$ of (\ref{eq:sigma}) by its associated limit quantity, as $N\rightarrow
\infty .$ We suppose from now on that the normalization condition 
\begin{equation}  \label{eq:rencN}
m! \left\Vert c_{N}\right\Vert _{2}^{2} = 1 , \mbox{ where } \left\Vert
c_{N}\right\Vert _{2}^{2}:=\sum_{\alpha \in \Lambda _{N,m}}c_{N}^{2}(\alpha
)
\end{equation}
is satisfied.
\newline
\textbf{Step 3.a (Law of large numbers)} Assume that $Y^{i},i=1,...,N,$ are
i.i.d. real valued random variables. We show that 
\begin{equation*}
\E\left( \sum_{\alpha \in \Lambda _{N,m}}c_{N}^{2}(\alpha
)(\prod_{i=1}^{m}f^{2}(Y^{\alpha _{i}})-\prod_{i=1}^{m}\E(f^{2}(Y^{\alpha
_{i}})))\right)^{2}\leq C\widehat{\delta }^2_{N}(c_{N}),
\end{equation*}%
where $\widehat{\delta }_{N}(c_{N})$ is a variant of the \textquotedblleft
influence factor\textquotedblright\ defined in (\ref{inf3}) below and where $ C$ depends on $m$ and on 
the supremum norm of $f.$ Since $%
Y^{i},i=1,...,N,$ are identically distributed, we have that $%
\E(f^{2}(Y^{\alpha _{i}}))=\E(f^{2}(Y^{1})),$ such that the above control
actually reads as 
\begin{equation*}
\E\left( \sigma_N^{2}(Y)-(\E f^{2}(Y^{1}))^{m}\right) ^{2}\leq C \widehat{%
\delta }^2_{N}(c_{N}).
\end{equation*}

At this stage we therefore have the following result. If $\overline{X}%
_{t_{k}}^{n,i},i=1,...,N, $ are i.i.d. for some $k, $ then 
\begin{equation*}
\E\left( (\overline{X}_{t_{k+1}}^{n,i} - \overline{X}_{t_{k}}^{n,i}- (\E f^{2}(%
\overline{X}_{t_{k}}^{n,i}))^{m/2}(B_{t_{k+1}}-B_{t_{k}})\right)^{2}\leq C %
\widehat{\delta }^2_{N}(c_{N}) \frac1n .
\end{equation*}%
However, the coordinates $\overline{X}_{t_{k}}^{n,i},i=1,...,N,$ are clearly
not i.i.d. for $k \geq 1, $ even when starting from i.i.d. initial
coordinates at time $t_0.$ Therefore, we have to adapt the above arguments,
and here is where the \textit{conditional propagation of chaos} property
comes in.

\noindent \textbf{Step 3.b (Conditional propagation of chaos) } We now introduce $b$ and $a$ again. We 
apply the usual crucial trick in the propagation of chaos theory by defining an
alternative chain 
\begin{multline}
\widehat{X}_{t_{k+1}}^{n,i}=\widehat{X}_{t_{k}}^{n,i}+\frac{1}{n} b(\widehat X_{t_k}^{n,i},(\E_{B}%
\varphi (\widehat{X}_{t_{k}}^{n,i}))^{\bar m})
+a(\widehat X_{t_k}^{n,i},(\E_{B}\varphi (\widehat{X}_{t_{k}}^{n,i}))^{\bar m})\frac{%
1}{\sqrt{n}}\Delta _{k}^{i}\\
+(\E_{B}f^{2}(\widehat{X}%
_{t_{k}}^{n,i}))^{m/2}(B_{t_{k+1}}-B_{t_{k}}),0\leq k<n.  \label{I4}
\end{multline}%
Notice that if $\widehat{X}_{{0}}^{n,i},i=1,...,N,$ are i.i.d., then for
every $k\geq 1 ,$ the $\widehat{X}_{t_{k}}^{n,i},i=1,...,N,$ are also
i.i.d., conditionally on $B.$ So we may use conditional independence and,
using the law of large numbers for polynomial interactions, we obtain 
\begin{eqnarray*}
\E(\sigma_N ^{2}(\widehat{X}_{t_{k}}^{n})-(\E_{B}f^{2}(\widehat{X}%
_{t_{k}}^{n,1}))^{m})^{2} &\leq &C \widehat{\delta }^2_{N}(c_{N}) , \\
\E(\mathbb{T} (\widehat{X}_{t_{k}}^{n})-(\E_{B}\varphi (%
\widehat{X}_{t_{k}}^{n,1}))^{\bar m})^{2} &\leq &C \widehat{\delta }^2_{N}(\bar c_{N})
\end{eqnarray*}%
(we recall that $\E_{B}$ denotes conditional expectation with respect to $%
\sigma \{B_{s},s\geq 0\}$).

Let us summarize what we have achieved so far. If $\widehat{X}_{0}^{n}=%
\overline{X}_{0}^{n}=X_{0}^{n}$ is composed of i.i.d. components, then, for 
a constant depending on $m$ and on the model parameters,  
\begin{multline*}
\mathcal{W}_{2}^{2}(\widehat{X}_{t_{k}}^{n},X_{t_{k}}^{n})\leq 2 \mathcal{W}%
_{2}^{2}(X_{t_{k}}^{n},\overline{X}_{t_{k}}^{n})+2 \mathcal{W}_{2}^{2}(%
\widehat{X}_{t_{k}}^{n},\overline{X}_{t_{k}}^{n})\\
\leq Ct_{k} [ \delta_N^p ( c_N) + \hat \delta_N^2 ( c_N) + 
\hat \delta_N^2 ( \bar c_N) +  \varepsilon_N ( c_N)] .
\end{multline*}
Since we supposed that $ \delta_N ( c_N), \delta_N ( \bar c_N) \le C N^{ - 1 / 2}, $ the above expression is of order 
$ N^{ - p /2 } +  \varepsilon_N ( c_N).$ 

\noindent \textbf{Step 4 } For every $i\in \{1,\ldots ,N\},$ the chain $%
\widehat{X}_{t_{k}}^{n,i},0\leq k\leq n,$ is the Euler scheme associated to
the stochastic equation 
\begin{eqnarray*}
\mathbb{X}_{t}^{i} &=&X_{0}^{i}+\int_{0}^{t}b(\mathbb{X}_{s}^{i},(\E_{B}(%
\varphi (\mathbb{X}_{s}^{i})))^{\bar m})ds \\
&&+\int_{0}^{t}a(\mathbb{X}_{s}^{i},(\E_{B}(\varphi (\mathbb{X}%
_{s}^{i})))^{\bar m})dW_{s}^{i}+\int_{0}^{t}(\E_{B}(f^{2}(\mathbb{X}%
_{s}^{i})))^{m/2}dB_{s},
\end{eqnarray*}%
Using the stochastic sewing lemma (see e.g. Lemma 2.3 in \cite{EvaVlad}), we
prove that this equation has a unique solution and that 
\begin{equation*}
\mathcal{W}_{2}^{2}(\widehat{X}_{1}^{n},\mathbb{X}_{1})\leq \frac{C}{%
{n}} .
\end{equation*}%
This concludes the sketch of the main ideas of the proof of Theorem \ref%
{theo:main0}.

We finish this section with a significant example. 
\begin{example}\label{ex:intro}
We consider $m=2$ and we take the
coefficients 
\begin{equation*}
c_{N}(i,j)=\kappa_{N}\frac{1}{2\sqrt{N\log N}}\times \frac{1}{\sqrt{\left\vert
j-i\right\vert }}\mathbf{1}_{\{i\neq j\}},
\end{equation*}%
where $\kappa_{N}\rightarrow 1$ as $N\rightarrow \infty $ is such that
the normalization relation
(\ref{eq:rencN}) is satisfied. Then the influence factor satisfies $\delta_N(c_N)\leq C/\sqrt{N},$ and the contraction term appearing in (\ref{M1}) satisfies $\left\Vert c_N\otimes_1 c_N\right\Vert_{1,2}^2\leq C\varepsilon _{N}(c_{N})$ with $\varepsilon _{N}(c_{N})=(\log N)^{-2}.$
Moreover, for $\theta >\frac{1}{2},$ we take 
\begin{equation*}
c_{N}(i,j)=\kappa_{N}\frac{1}{2\sqrt{J(2\theta)N}}\times \left\vert j-i\right\vert
^{-\theta}\mathbf{1}_{\{i\neq j\}},
\end{equation*}%
with $J(2\theta)=\sum_{k=1}^{\infty }k^{-2\theta}.$ Then $\delta_N(c_N)\leq C/\sqrt{N},$ and the same holds with 
\begin{equation*}
\varepsilon _{N}(c_{N})=\left\{ 
\begin{array}{cc}
\frac{1}{N^{4\theta-2}}, & \mbox{ if }\frac{1}{2}<\theta<\frac{3}{4} \\ 
\frac{\log N}{N}, & \mbox{ if }\theta=\frac{3}{4} \\ 
\frac{1}{N}, & \mbox{ if }\theta>\frac{3}{4}%
\end{array}%
\right\} .
\end{equation*}%
In particular, the right hand side of (\ref{M1}) is bounded by $C_p(r_N+\frac{1}{n}),$ with $r_N=(\log N)^{-2}$ for $\theta=\frac{1}{2},$ $r_N=N^{-p/2}+N^{-(4\theta-2)}$ for $\frac{1}{2}<\theta<\frac{3}{4}$ and $r_N=N^{-p/2}$ for $\theta\geq\frac{3}{4}$ (see Theorem \ref{theo:example} below).

We mention that one may also consider the case $%
m\geq 3,$ but then concrete computations become heavier.
\end{example}

\subsection*{Bibliographic remarks}
McKean-Vlasov equations with common noise and the conditional propagation of chaos
property have been studied by many authors, see for instance Kurtz and Xiong \cite{KX99},
Coghi and Flandoli \cite{CF16}, Carmona, Delarue and Lacker \cite{CDL16} and the monograph
\cite{CD18}. In these works, the common noise is already present in the finite particle
system. In contrast to this, in the present paper, as in \cite{ELL21}, \cite{ELL23} and \cite{EvaVlad}, the
common noise is not present in the finite system: it emerges in the limit, as the
cumulative effect of many small random inputs to which all particles are equally
exposed. In \cite{ELL21} this has been shown for mean field systems of interacting
neurons in a diffusive scaling, where the Brownian common noise is created by the
CLT. Strong error bounds for this model have then been given in \cite{ELL23}, and a
general quantitative framework has been developed in \cite{EvaVlad}. In all these works, the
bombardment term is a linear sum. In contrast to this, in the present work it is a multilinear homogeneous sum.

Limit theorems for multilinear forms and homogeneous sums go back to Rotar'
\cite{R75}, \cite{R79} and de Jong \cite{dJ89}, \cite{dJ90}. The invariance principle for multilinear
polynomials with low influences is due to Mossel, O'Donnell and Oleszkiewicz \cite{MDO},
see also Nourdin, Peccati and Reinert \cite{NPR} for the link with Wiener chaos, and Bally
and Caramellino \cite{BC} for estimates in total variation distance. For multiple
stochastic integrals, the fourth moment theorem of Nualart and Peccati \cite{NuPe} implies
that a sequence of multiple integrals of fixed order converges to a normal law if
and only if the associated contractions vanish. Quantitative versions of this result
are given in \cite{NP} and \cite{LNP}. We use these results at each time step, with coefficients
which depend on the configuration of the system.

Finally, our system is not exchangeable, because the weights $c_N(\alpha)$ and $ \bar c_N(\alpha)$ are not required to be invariant under a
relabelling of the particles. However, every particle feels the rest of the system only through two global quantities: the interaction term $ \mathbb{T}(x)$ and the collective perturbation $S_k(x).$
Both are the same for all particles. In this sense, the interaction is still of a “global” type, even though it is not a function of the empirical measure. This is why the spatial structure of the weights disappears in the limit and the limit equation (\ref{eq:SDE}) is the same for all particles. This differs from the non-exchangeable systems studied, for instance, by Jabin, Poyato and Soler \cite{JPS25}, where each particle $i$ interacts with the others through its own weights $w_{ij},$ and where the limit keeps track of the label of the particle. Extending our results to particle-dependent interactions, such as $ \mathbb{T}_i(x) = \sum_{ \alpha : i \in \alpha} \bar c_N^2 (\alpha) \prod_{j=1}^{\bar m} \varphi (x^{ \alpha_j} )  ,$ is a natural question that we leave for future work.

\subsection*{Our contribution}
The general strategy of this paper is to couple
the particle system with an Euler scheme of the limit equation, controlling the
error by a stochastic sewing argument. In \cite{EvaVlad}, the common noise is created by a
linear sum of small independent inputs, and the Gaussian limit is obtained by a
quantitative CLT. The new point here is that the bombardment term $S_k$ is a
homogeneous sum of order $m$, that is, a polynomial of the inputs, which is
evaluated along the configuration of the system. To deal with this, we combine
three ingredients. First, we use the invariance principle of \cite{MDO}, in the form given
in \cite{BC} and \cite{NPR}, and the CLT for multiple stochastic integrals of \cite{LNP} and \cite{NP}, with coefficients
which depend on the current configuration $x$. Second, since the interactions are
not of mean field type and the system is in general not exchangeable, we prove a
law of large numbers for polynomial interactions (Proposition \ref{prop:LLN}). We stress that this step is crucial 
since it allows us to switch to the limit exchangeable system. Third, we
iterate these one-step estimates along the $n$ time steps, following \cite{EvaVlad}, without
losing a factor $n^{-1/2}$ at each step. We do not claim that our assumptions or our
rates are optimal. In particular, the restriction $p < 1/3$ in (\ref{M1}) comes from the
passage from $d_3$ to $\mathcal{W}_2$, and our method does not cover the case
$0 < \theta < 1/2$ in Example \ref{ex:intro}.

\subsection*{Notation}

Throughout this paper, $C$ and $K$ denote arbitrary positive constants whose
values can change from line to line in an equation. We write $C_\theta$ and $%
K_\theta$ if the constants depend on some parameter $\theta.$

The following objects will be frequently used throughout the paper, for some fixed $ m, N \geq 1. $  We denote \\
$-$  $\Lambda _{N, m}$ the set of multi-indices $\alpha
=(\alpha _{1},...,\alpha _{ m})\in \{1,...,N\}^{m};$ \\
$-$ $\Lambda _{N, m}^{i}$\ the set of the multi-indices which
contain $i$ (for some $ 1 \le i \le N $ which is fixed);\\
$-$ $%
c_{N}(\alpha )\in {\mathbb{R}},\alpha \in \Lambda _{N,m},$ a collection of non-negative coefficients which are
symmetric (invariant under permutation of its arguments) and null
on the diagonals;\\
$-$ the {\it influence factor} $\delta _{N}(c)$  given by 
\begin{equation}
\delta _{N}^{2}(c_{N})=\max_{i=1,...,N}\sum_{\alpha \in \Lambda
_{N, m}^{i}} c_{N}^{2}(\alpha ) ;  \label{inf2}
\end{equation}%
$-$ the {\it modified influence factor} \begin{equation}
\widehat{\delta }_{N}^{2}( c_{N})=\max_{k=1,\ldots ,  m}\sum_{i=1}^{N}\left(
\sum_{\alpha \in \Lambda _{N, m},\alpha _{k}=i} c_{N}^{2}(\alpha )\right)
^{2}\leq N\delta _{N}^{4}(c_{N});  \label{inf3}
\end{equation}
$-$ the \textit{contraction factors} (see \cite%
{NuPe}) given by 
\begin{equation}
\left\Vert c_{N}\otimes _{q}c_{N}\right\Vert _{m-q,2}^{2}=\sum_{\alpha
,\beta \in \Lambda _{N,m-q}}\left( \sum_{\gamma \in \Lambda
_{N,q}}c_{N}(\alpha ,\gamma )c_{N}(\beta ,\gamma )\right) ^{2},  \label{inf4}
\end{equation}
for all $1\leq q\leq m-1.$ 

\section{Interacting particle systems with polynomial interactions}\label{sec:2}
In this section we present the general dynamic defining the one step Euler scheme for the McKean-Vlasov type
equation with polynomial interactions we are interested in. This dynamic is a variant of the model 
considered in \cite{EvaVlad}. To begin with, 
we introduce the general framework we are working in. We work with a fixed number of particles $N.$

\subsection{Framework}\label{sec:frameworkgeneral} To represent the processes we are interested in
we fix a generic probability space $(\Omega ,\mathcal{A},P).$ All random variables 
we shall consider are defined on this space. We fix two times $ 0 \le s < t .$ Our assumptions will be stated for these fixed
times. 

We consider five sub-$\sigma-$algebras of $ \mathcal A:$ $\mathcal F_s^\Delta \subset \mathcal F_s \subset  \mathcal F_t $ and  $\mathcal{G}^1_{s,t},  \mathcal G^2_{s, t } \subset \mathcal F_t .$
We also consider a standard Gaussian random variable $ \Delta_{s, t } $ which is  $ \mathcal F_t-$measurable. 
We suppose that $\mathcal{G}^i_{s,t}, i= 1, 2, \mathcal F_s$ and $ \sigma ( \Delta_{s, t } ) $ are mutually independent.
The significance of $ \mathcal F_s^\Delta$ will become clear when we iterate the one-step scheme over finite partitions, see (\ref{eq:Ftdelta}) below.

The two $\sigma-$algebras  $\mathcal{G}_{s,t}^{1}$ and $\mathcal{G}_{s,t}^{2}$ have to 
be thought of representing two sources of {\it innovation} within the period $(s,t),$
which are independent of the past $\mathcal{F}_{s}.$ The first source of innovation, which is $\mathcal{G}_{s,t}^{1}-
$adapted, will be related to an autonomous evolution of the different
particles, while the second source of innovation,
which is $\mathcal{G}_{s,t}^{2}-$adapted, will be a perturbation term to
which each particle is equally exposed. 

For any $ n \geq 1, $ we denote by $L_{s}^{p} (\R^n ) $ the space of the $\mathcal{F}_{s}-$measurable ${\mathbb{R}}^{n}-$valued 
random variables with coordinates which all have a finite moment of order $p,$ and more generally 
by $ L^p ( \R^n )$ the space of all  ${\mathbb{R}}^{n}-$valued 
random variables with coordinates which all have a finite moment of order $p.$
We denote by $\E_{s}(X)$
the conditional expectation of $X$ with respect to $\mathcal{F}_{s},$ and we
put $\widetilde{\E}_{s}(X)=X-\E_{s}(X).$ The main difference with respect to
the presentation in \cite{EvaVlad}  is that there one
works with a filtration $\mathcal{F}_{t},t\geq 0,$ while here $s$ and $t$ are
fixed. 

We introduce now the core objects of our dynamics.

$\clubsuit $ \qquad {\bf Autonomous innovation.} 
We consider an i.i.d. collection $H_{s,t}^{i}:\Omega \times {\mathbb{R}}\times {\mathbb{R}}%
\rightarrow {\mathbb{R}},i=1,...,N,$ of random fields such that each $H_{s,t}^{i}$ is $\mathcal{G}^1_{s,t}
\otimes \mathcal{B} ( {\mathbb{R}} \times {\mathbb{R}}) - \mathcal{B} ( {
\mathbb{R}} )-$measurable, and such that 
\begin{equation}  \label{eq:H}
H_{s,t}^{i}( X,a)\in L_{t}^{2}({\mathbb{R}}) \mbox{ whenever  } ~X\in
L_{s}^{2}({\mathbb{R}}) \mbox{ and  } a\in L_{s}^{2}({\mathbb{R}}) .
\end{equation}%
We put $H_{s,t}=(H_{s,t}^{1},...,H_{s,t}^{N})\in {\mathbb{R}}^{ N}.$

We impose the following assumptions on $ H_{s, t}$ which are natural and 
straightforward to check in most examples. 

\begin{assumption}\label{ass:H}

\begin{enumerate}
\item If $X^{1},...,X^{N} \in L_{s}^{2}({\mathbb{R}})$ are i.i.d. and if 
$a\in {\mathbb{R}}$ is deterministic, then $H_{s,t}^{i}(X^{i},a),$ $%
i=1,...,N,$ are i.i.d.

\item We assume that there exists a constant $C_{stb}$\ such that, almost
surely, for every $x,y\in {\mathbb{R}},a,b\in {\mathbb{R}},$ 
\begin{equation}
\left\vert \E_s\left(H_{s,t}^{i}(x,a)-H_{s,t}^{i}(y,b)\right)\right\vert \leq
C_{stb}(t-s)(\left\vert x-y\right\vert +\left\vert a-b\right\vert ).
\label{pf1}
\end{equation}%
Moreover, for every $X^i,Y^i\in L_{s}^{2}({\mathbb{R}}),A,B\in L_{s}^{2}(%
{\mathbb{R}}),$%
\begin{equation}
\E\left\vert \widetilde{\E}_{s}(H_{s,t}^{i}(X^i,A)-H_{s,t}^{i}(Y^i,B))\right%
\vert ^{2}\leq C^2_{stb}(t-s)(\E\left\vert X^i-Y^i\right\vert
^{2}+\E\left\vert A-B\right\vert ^{2}) .  \label{pf1'}
\end{equation}
\end{enumerate}
\end{assumption}

Here is a typical example which we quote from our previous paper \cite{EvaVlad}. 
\begin{example}
Let $N^i, 1 \le i \le N, $ be independent Poisson random measures on ${%
\mathbb{R}}_+ \times {\mathbb{R}} $ having intensity $dt \nu ( dv) $ where 
$\nu $ is a $\sigma-$finite measure having a finite first and second moment. Let
moreover $\beta^i, 1 \le i \le N, $ be independent $1-$dimensional Brownian
motions, independent of the Poisson random measures. Let $\psi : {\mathbb{R}}%
 \times {\mathbb{R}} \to {\mathbb{R}} $ and $\kappa : {\mathbb{R}}
\to {\mathbb{R}} $ be Lipschitz continuous and consider 
\begin{equation*}
H_{s,t}^{i}(x,a)=\psi (x,a)(t-s)+\int_{] s, t ] } \int_{{\mathbb{R}}
}vN^{i}(dr,dv)+\kappa (x)(\beta_{t}^{i}-\beta_{s}^{i}).
\end{equation*}
\end{example}

\subsection*{The interactions}
To describe the autonomous interactions in the system, we fix $\varphi :\R\rightarrow \R ,$
a bounded Lipschitz continuous function with Lipschitz constant $%
L_{\varphi },$ and we fix some $\bar m \geq 1 ,$ the order of the interactions. 
Mean field interactions correspond to the choice of $ \bar m = 1 , $ yielding an interaction of type 
$\frac{1}{N}\sum_{i=1}^{N}\varphi (X^{i})$ 
appearing in standard McKean-Vlasov equations. 
We shall consider polynomial interactions of order $ \bar m$  which 
are a generalization of
these mean field interactions. Recall that 
 $\Lambda _{N,\bar m}$ denotes the set of multi-indices $\alpha
=(\alpha _{1},...,\alpha _{\bar m})\in \{1,...,N\}^{\bar m}.$ For any vector $%
x=(x^{1},...,x^{N})\in {\mathbb{R}}^{N}$ and any $ \alpha \in \Lambda_{N, \bar m}, $  
we then denote $x^{\alpha
}=x^{\alpha _{1}}\times ...\times x^{\alpha _{\bar m}}.$

In what follows we consider a collection of non-negative symmetric coefficients $\bar c_{N}(\alpha
),\alpha \in \Lambda _{N,\bar m},$ which are null on the diagonal. 
We suppose that the normalization condition
\begin{equation}\label{eq:barcnorm}
\left\Vert 
\overline{c}_{N}\right\Vert^2_{2}= \sum_{\alpha \in \Lambda
_{N,\bar m}}\bar c_{N}^{2}(\alpha ) = 1 
\end{equation} 
holds. We then define%
\begin{equation}
\mathbb{T} (x)=\sum_{\alpha \in \Lambda _{N,\overline{m}}}\overline{c}%
_{N}^{2}(\alpha )\prod_{j=1}^{\overline{m}}\varphi (x^{\alpha _{j}}) 
\label{McKV6'}
\end{equation}
and we put 
\begin{equation}
h_{s,t}^{i}(x)=H_{s,t}^{i}(x^{i},\mathbb{T} (x)).  \label{McKV7''}
\end{equation}
$ h_{s, t }^i ( x) $ is the autonomous evolution of particle $i, $ between times $s$ and $t, $ if
at time $s, $ the initial configuration equals $x.$ 

We have in mind two significant examples.

\begin{example}

\textbf{1. Mean field interaction.} We take $\overline{m}=1$ and $\overline{c}%
_{N}(i)=\frac{1}{\sqrt{N}}$ such that $\mathbb{T} (x)=\frac{1}{N}%
\sum_{j=1}^{N}\varphi (x^{j}).$

\textbf{2. Polynomial interaction.} We take $\overline{m}=2$%
\textbf{\ }and $\overline{c}_{N}(i,j)= \kappa_N \frac{1}{\sqrt{2 N \log N}}\times \frac{1}{%
\sqrt{\left\vert i-j\right\vert }}\times \mathbf{1}_{i\neq j},$ where $ \kappa_N $ is such that $ \| \bar c_N\|_2^2 = 1 .$ 
In this case
\begin{equation*}
\mathbb{T} (x)=\frac{\kappa_N^2 }{2 N\log N}\sum_{i, j=1, i \neq j }^{N}\frac{\varphi
(x^{i})\varphi (x^{j})}{\left\vert i-j\right\vert }.
\end{equation*}
\end{example}

Suppose now that $X^i \in L^2_s ( {\mathbb{R}}) $
for any $1 \le i \le N, $ and that the initial configuration is random, given by $X = ( X^1, \ldots, X^N).$
In this case, we also introduce $\widehat{h}_{s,t}^{i}:L_{s}^{2}({\mathbb{R%
}}^{ N})\rightarrow L_{t}^{2}({\mathbb{R%
}}^{ N}),i=1,...,N,$ by 
\begin{equation}
\widehat{h}_{s,t}^{i}(X)=H_{s,t}^{i}(X^{i},(\E(\varphi (X^{i})|\mathcal{F}^\Delta
_{s})^{\bar m } )) .  \label{pf3}
\end{equation}%
We write $h_{s, t}(X)=(h_{s, t}^{1}(X),...,h_{s,t}^{N}(X))\in {\mathbb{R}}^{ N},$ and $%
\widehat{h}_{s, t }(X)=(\widehat{h}_{s, t }^{1}(X),...,\widehat{h}_{s, t} ^{N}(X)).$ Notice that if $%
X = (X^1, \ldots, X^N) $ is exchangeable, then  $\widehat{h}_{s, t} 
( X )$ is exchangeable, too. This is in contrast to $ h_{s, t} ( X) $ which is in general not exchangeable. 

$\clubsuit $ \qquad \textbf{Perturbations.} In our dynamics, each particle is equally exposed to some 
random perturbation that we are going to describe now. We
consider an i.i.d. collection of random fields 
\begin{equation}  \label{pf4}
v_{s,t}^{i}(x^i) = v_{s,t}^{i}(x^i,\omega ) \in {\mathbb{R}},\;
x^i \in \R, \omega \in \Omega , 1 \le i \le N , 
\end{equation}
satisfying the following assumption.

\begin{assumption}
\label{ass:perturbation} i) For any $ 1 \le i \le N, $ we have that 
\begin{equation*}
\Omega \times {\mathbb{R}}\ni (\omega ,x^i )\mapsto v_{s, t}^i
(x^i, \omega )\in {\mathbb{R}}
\end{equation*}%
is $\mathcal{G}_{s,t}^{2}\otimes \mathcal{B}({\mathbb{R}})-%
\mathcal{B}({\mathbb{R}})-$ measurable.

ii) For any fixed $x = (x^1, \ldots, x^N ) \in {\mathbb{R}}^N,$ the random variables $%
v_{s, t}^i ( x^i, \omega) , 1 \le i \le N ,$ are independent. Moreover, they are all
centered and square integrable. 

iii) If $ X^1, \ldots, X^N \in L^2_s ( \R ) $ are i.i.d., then the $v^i_{s, t } (X^i ), 1 \le i \le N ,$ are i.i.d. as well. 

iv) There exists a Lipschitz continuous function $f : \R \to \R_+$ such that for all $1 \le i \le N, $ for all $ x^i \in \R,$
\begin{equation*}
\E(v_{s, t} ^{i}(x^{i}))=0,\quad \E(v_{s, t} ^{i}(x^{i})^{2})=f^{2}(x^{i}) . 
\end{equation*}%
Finally we also suppose that for all $x^i \in \R,$ 
\begin{equation*}
0<\underline{f}\leq f(x^i )\leq \overline{f}.
\end{equation*}
\end{assumption}
Since $f$ is Lipschitz and bounded, $f^2 $ is Lipschitz continuous and bounded as well. In what follows, we shall write $L_{f^2 } $ for the associated 
Lipschitz norm, that is, $ |f^2 (u) - f^2 (v) | \le L_{f^2 } | u-v | $ for all $ u, v \in \R.$ 

We stress that since the $ v^i_{s, t } ( x^i , \omega ) $ are $\mathcal{G}_{s,t}^{2}\otimes 
\emph{B}({\mathbb{R}})$\ measurable, they are independent of the past $%
\mathcal{F}_{s}$ and independent of $ \mathcal{G}_{s,t }^1.$ In particular, they are independent of the variables $h^i_{s,t} (x) $
introduced above.

For a fixed $q\geq 3$ we consider the following moment assumption:

\begin{assumption}\label{ass:2.10first}
We assume that there exists $\gamma_{s, t} (q)<\infty $ such that
for every $1 \le i \le N $ and every $x^i \in \mathbb{R},$%
\begin{equation}
E(\left\vert v_{s, t }^{i}(x^i)\right\vert ^{q})\leq \gamma_{s, t} (q).  \label{CLT2first}
\end{equation}
\end{assumption}

To illustrate the objects we are interested in, let us quote the following example from our previous paper \cite{EvaVlad}. 

\begin{example}
Recall that $ s < t $ are fixed. Let $M_{s, t } ^{i}(dz, du),1\leq i\leq N,$ be i.i.d. Poisson random measures on ${%
\mathbb{R}}_{+}\times {\mathbb{R}}$ having
intensity $dz\pi (du)$ such that $\pi $ is centered with 
\begin{equation*}
\int_{{\mathbb{R}}} u^2 \pi (du)=1 ,\end{equation*}%
having finite moments of order $q.$ 
Putting 
\begin{equation*}
v_{s,t}^{i}(\omega ,x^{i})=\int_{{\mathbb{R}}_{+}}\int_{{\mathbb{%
R}}}u1_{\{z\leq f^2  (x^{i})\}}M_{s, t }^{i}(dz,du) ,
\end{equation*}%
it is easy to see that the above assumptions are satisfied.
\end{example}

\subsection*{Defining the polynomial interactions}

As above, we fix an order of interactions $ m \geq 2$ and we consider another set of coefficients $c_{N}(\alpha ),\alpha \in \Lambda
_{N,m},$  satisfying the
normalization condition%
$$
m!\left\Vert c_{N}\right\Vert _{2}^{2}=1  
$$
and
\begin{equation}
\delta _{N}(c_{N})\leq \frac{C}{\sqrt{N}} \label{McKV8'}
\end{equation}%
(recall (\ref{inf2})). 

For any $\alpha =(\alpha _{1},...,\alpha _{m})\in \Lambda _{N,m}$ and $%
x=(x^{1},...,x^{N})\in {\mathbb{R}}^{N},$ we put%
\begin{equation*}
v_{s, t} ^{\alpha }(x)=v_{s, t} ^{\alpha _{1}}(x^{\alpha _{1}})\times ...\times v_{s, t} ^{\alpha
_{m}}(x^{\alpha _{m}}),
\end{equation*}%
and we associate to the vector $v_{s, t} (x)=(v_{s, t}^{1}(x^{1}),...,v_{s, t} ^{N}(x^{N}))$ the random variable
\begin{equation}
S_{s,t}^{N}(c_{N},x)=\sqrt{t-s}  \sum_{\alpha \in \Lambda
_{N,m}}c_{N}(\alpha )v_{s,t}^{\alpha }(x).  \label{McKV8'''}
\end{equation}%

We also fix some standard normal random variable $\Delta _{s,t},$ independent of $%
\mathcal{F}_{s},$ and we put 
\begin{equation}
\overline{S}_{s,t}^{N}(c_{N},x)= \sigma _{N}(x)\sqrt{t-s}\Delta _{s,t} ,
\label{McKV8}
\end{equation}%
where 
\begin{equation}
\sigma _{N}^{2}(x)=\E\left( \left(\sum_{\alpha \in \Lambda
_{N,m}}c_{N}(\alpha )v_{s,t}^{\alpha }(x)\right)^2\right) =m! \sum_{\alpha \in \Lambda
_{N,m}}c_{N}^{2}(\alpha )\prod_{i=1}^{m}f^2 (x^{\alpha _{i}}). \label{McKV9}
\end{equation}

Finally, if $X = ( X^{1},\ldots ,X^{N}) \in L^2_s( {\mathbb{R}}^{ N} ),$ we also put
\begin{equation}
\widehat{S}_{s,t}^{N}(c_N, X)= (\E (f^2 ( X^1) | \mathcal{F}^\Delta_s
))^{m/2}  \sqrt{t-s} \Delta_{s, t} .   \label{Avibis}
\end{equation}

\subsection*{The one-step dynamic}
We are now able to define the one-step dynamic of each particle which is given by the application
\begin{equation}
\mathcal{X}_{s,t}^{i}(X)=X^{i}+h_{s,t}^{i}(X)+S_{s,t}^{N}(c_{N},X) ,
\label{McKV11}
\end{equation}%
for any $ 1 \le i \le N$ and for any $ X \in L^2_s ( \R^N)  .$  
Moreover, we have the two approximating schemes
\begin{equation}
\overline{\mathcal{X}}_{s,t}^{i}(X)=X^{i}+h_{s,t}^{i}(X)+\overline{S}%
_{s,t}^{N}(c_{N},X), \label{McKV12}
\end{equation}
and 
\begin{equation}
\widehat{\mathcal{X}}_{s,t}^{i}(X)=X^{i}+\widehat h_{s,t}^{i}(X)+\widehat{S}%
_{s,t}^{N}(c_{N},X), \label{McKV12hat}
\end{equation}
which are defined 
for any $ 1 \le i \le N$ and for any $ X \in L^2_s ( \R^N)  .$

\subsection{Iterating the one-step scheme over finite partitions}
We consider now the time grid $\pi
=\{t_{0}<t_{1}<...<t_{k}<...\},$ starting from $t_0 = 0,$ and we denote $\gamma _{k}=t_{k+1}-t_{k}.$
Let us specify the framework associated to this time grid. We assume that 
we have a filtration $(\mathcal{F}_{t_{k}})_{ k \geq 0}.$ 
Moreover we assume that for all $ k \geq 0, $ we have some $\sigma -$algebras $%
\mathcal{G}_{t_{k},t_{k+1}}^{i}\subset \mathcal{F}_{t_{k+1}},i=1,2,$ such that $ \mathcal{ F}_{t_k}, 
\mathcal{F}_{t_{k+1}} $ and $ \mathcal{ G}^i_{t_k, t_{k+1} }, i=1,2, $ satisfy the assumptions of Section \ref{sec:frameworkgeneral}
above, with $s := t_k$ and $ t := t_{k+1}.$ We also suppose that  Assumptions \ref{ass:H}, \ref{ass:perturbation} and \ref{ass:2.10first} are satisfied for every $k \geq 0, $ 
choosing again $s := t_k$ and $ t:= t_{k+1}.$ We put \begin{equation}\label{eq:gammaTq}
 \gamma_T ( q) := \sup_{ k : t_{k+1} \le T } \gamma_{t_k, t_{k+1} } ( q) .
\end{equation}

As before, for any $ X \in L^1 ( \R), $ we
denote $\E_{t_{k}}(X)$ the conditional expectation of $X$ with respect to $%
\mathcal{F}_{t_{k}}$ and $\widetilde{\E}_{t_{k}}(X)=X-\E_{t_{k}}(X).$

We also consider an i.i.d. sequence of standard normal  random
variables $\Delta _{k}=\Delta _{t_{k},t_{k+1}},k\geq 0, $ such that $\Delta_k$ 
is $ \mathcal{F}_{t_{k+1}}-$measurable and independent of $ \mathcal{F}_{t_k}.$
We denote
\begin{equation}\label{eq:Ftdelta}
\mathcal{F}_{t_{k}}^{\Delta }=\sigma (\Delta _{0},...,\Delta _{k-1})
\end{equation} 
and $%
\mathcal{F}^{\Delta }=\cup _{k\geq 0}\mathcal{F}_{t_{k}}^{\Delta }.$
Notice that $ \mathcal F_{t_k}^\Delta \subset \mathcal F_{t_k}.$
We also
denote by $\E_{t_{k}}^{\Delta }(X)$ the conditional expectation with respect
to $\mathcal{F}_{t_{k}}^{\Delta }$ (which is of course different from $%
\E_{t_{k}}(X)).$ 

Then we consider the chains%
\begin{eqnarray}
X^i_{t_{k+1}} &=&\mathcal{X}^i%
_{t_{k},t_{k+1}}(X_{t_{k}})=X_{t_{k}}^{i}+h_{t_{k},t_{k+1}}^{i}(X_{t_{k}})+S_{t_{k},t_{k+1}}^{N}(c_{N},X_{t_k}),
\notag \\
\overline{X}^i_{t_{k+1}} &=&\overline{\mathcal{X}}^i_{t_{k},t_{k+1}}(\overline{X}%
_{t_{k}})=\overline{X}_{t_{k}}^{i}+h_{t_{k},t_{k+1}}^{i}(\overline{X}%
_{t_{k}})+\sigma _{N}(\overline{X}_{t_{k}})\sqrt{t_{k+1}-t_{k}}\Delta
_{t_{k},t_{k+1}}, \notag \\
\widehat{X}^i_{t_{k+1}}&=&\widehat{\mathcal{X}}^i_{t_{k},t_{k+1}}(\widehat{X}%
_{t_{k}})=\widehat{X}_{t_{k}}^{i}+\widehat{h}_{t_k,t_{k+1}}^{i}(\widehat{X}%
_{t_{k}})\notag \\
&& \quad \quad \quad +(\E_{t_{k}}^{\Delta }f^{2}(\widehat{X}^1_{t_{k}}))^{m/2}\sqrt{t_{k+1}-t_{k}} \Delta
_{t_{k},t_{k+1}} .  \label{McKV13} 
\end{eqnarray}

The chain $ X_{t_k}, k \geq 0, $ describes our original particle system evolving in discrete time. 
The chain $ \overline X_{t_k}, k \geq 0, $ is an auxiliary chain, and $ \widehat X_{t_k}, k \geq 0, $ is the limit dynamic in discrete time.

Our main result now reads as follows. We fix some time horizon $ T > 0.$ 

\begin{theorem}
\label{theo:main} We suppose that there exists some constant $C$ such that 
\begin{equation}
\delta _{N}(\overline{c}_{N}) + \delta _{N}({c}_{N}) \leq \frac{C}{\sqrt{N}},  \label{McKV6}
\end{equation}
and we recall that $ m! \|c_N\|_2^2 = 1 $ and $ \| \bar c_N \|_2^2 = 1.$ 
Furthermore we impose Assumptions \ref{ass:H} and  %
\ref{ass:perturbation} for all $ (s, t) = (t_i, t_{i+1}), $ such that $ t_{i+1} \le T.$ 
Finally, we also suppose that  for some $ q \geq 3,$ Assumption  \ref{ass:2.10first} is verified for all $ (s, t) = (t_i, t_{i+1}), $ such that $ t_{i+1} \le T.$

Let $X_{0}^{1},...,X_{0}^{N} \in L^2_0 ( \R) $ 
be independent and identically distributed random variables, and suppose that $ X_0 = \overline X_0 = \widehat X_0.$  Then 
the following assertion holds. For every 
\begin{equation}
p<p(q):=\frac{2(q-2)}{6(q-2)+6},  \label{pq}
\end{equation}
there exists a constant $Q $ depending on the model parameters, on $ \gamma_T ( q),$ $ m, \bar m, $ $q, p $ and on the constant $C$ of (\ref{McKV6}), 
such that for every $ t_k \le T,$ 
\begin{equation}
\mathcal{W}_{2}^{2}(X_{t_k},\widehat{X}_{t_k})\leq
Q e^{QT^{2}}(N^{- p/2}  + \sum_{r=1}^{ m-1} \| c_N \otimes_r c_N\|_{m-r, 2 }^2) .  \label{pf10}
\end{equation}%
\end{theorem}

The proof of this theorem follows the ideas explained in the Introduction, and the remainder of the paper is devoted to it.

\subsection{From the Euler scheme to the limit equation}
We now come back to the particular framework of system (\ref{I1}) considered
in the introduction. As explained there, the discrete time chain $ (\widehat X_{t_k})_k$ appears as the Euler scheme 
approximation of the strong solution of equation (\ref{eq:SDE}). Strong existence and uniqueness for this equation have already been proved in 
\cite{EvaVlad}. The estimate of the error for the Euler approximation is based on the stochastic sewing lemma of L\^e \cite{Le}, rephrased for flows on the Wasserstein space
 by Alfonsi et al. \cite{ABC}. The arguments and the statement of this result are a slight variant of the one given in \cite{EvaVlad}, so in what follows we only recall them shortly. 

\begin{theorem}\label{theo:mainbis}
We assume that the hypotheses from the previous theorem are in force and we
take 
\begin{equation*}
H_{s,t}^{i}(x,u)=b(x,u)(t-s)+a(x,u)(W_{t}^{i}-W_{s}^{i}),
\end{equation*}
with $a$ and $b$ Lipschitz continuous and $\varphi $ Lipschitz continuous
and bounded. Then 
\begin{equation*}
\mathcal{W}_{2}^{2}(X_{1},\mathbb{X}_{1})\leq C\Big(N^{-p/2}+%
\sum_{r=1}^{m-1}\Vert c_{N}\otimes _{r}c_{N}\Vert _{m-r,2}^{2}+n^{-1}\Big) .
\end{equation*}
\end{theorem}

\begin{proof} 
We only have to compare the limit chain $\widehat X$ with the solution $\mathbb X$ of
(\ref{eq:SDE}). This is done exactly as in Section 6 of \cite{EvaVlad}, and we only explain
how our framework fits into the one of \cite{EvaVlad}. 
We put
$\Delta^i_k = \sqrt{n}\,(W^i_{t_{k+1}} - W^i_{t_k})$ and
$\Delta_k = \sqrt{n}\,(B_{t_{k+1}} - B_{t_k})$, where $W^i, i \ge 1,$ and $B$ are the
independent Brownian motions of (\ref{eq:SDE}). 

Since the future increments of $B$
are independent of the chain up to time $t_k$, $\E^\Delta_{t_k}$ coincides with the
conditional expectation $\E_B$ appearing in (\ref{eq:SDE}). Moreover, by construction,
$\widehat X^i_{t_{k+1}}$ only depends on $\widehat X^i_{t_k}$, on the conditional law of
$\widehat X^1_{t_k}$ given $\mathcal F^\Delta_{t_k}$, and on the increments of $W^i$ and $B$
on $[t_k, t_{k+1}]$ (see also Remark 3.8 of \cite{EvaVlad}). So, for each fixed $i$,
$\widehat X^i_{t_k}, 0 \le k \le n,$ is the Euler scheme of (\ref{eq:SDE}) along the grid
$\pi = \{ t_k = k/n, 0 \le k \le n \}$, and its one step scheme is the scheme which is denoted 
$\widehat{\mathcal X}_{s,t}$ in Section 3 of \cite{EvaVlad}, when taking $d = d_0 = 1,$ and choosing 
$$\widehat S_{s,t} (X) = ( \E ( f^2 ( X^1) | \mathcal F^B_s ))^{m/2} (B_t - B_s).$$
In other words, Assumption 3.5 of \cite{EvaVlad} holds with $f$ replaced by $f^2$ and with
$\kappa = m/2$. 

Moreover, since $\varphi$ is bounded and since
$\underline f^2 \le f^2 \le \overline f^2$, it follows that the maps $u \mapsto u^{\bar m},$ in restriction to
$[ - \|\varphi\|_\infty, \|\varphi\|_\infty ],$ and $u \mapsto u^{m/2},$ in restriction to 
$[\underline f^2, \overline f^2],$ are Lipschitz continuous. 

Hence the coefficients of
(\ref{eq:SDE}) satisfy the Lipschitz and boundedness conditions used in Section 6 of
\cite{EvaVlad}. In particular, the arguments of Lemma 6.2 and Proposition 6.3 of
\cite{EvaVlad} apply without change: $\widehat{\mathcal X}$ satisfies assumption
$\mathbf A(\beta)$ of Section 2 of \cite{EvaVlad} with $ \beta = 1.$  Moreover, its flow solution is the unique
solution of (\ref{eq:SDE}). In particular, (\ref{eq:SDE}) is well posed. Item {\bf C} of the
stochastic sewing lemma (Lemma 2.3 of \cite{EvaVlad}), applied to the partition $\pi$ of
$[0,1]$ with $|\pi| = 1/n$, then gives, for each $i$,
\[
\E | \mathbb X^i_1 - \widehat X^i_{1} |^2
\le C(\beta) \left(1 + E|X^i_0|^2\right) n^{-1 } .
\]
Since $\widehat X$ and $\mathbb X$ are constructed with the same $X_0$, $W^i$ and $B$, we
deduce, averaging over $i$, that
$\mathcal W_2^2 ( \widehat X_1, \mathbb X_1 ) \le E | \widehat X_1 - \mathbb X_1 |_N^2
\le C(\beta)(1 + E|X^1_0|^2) n^{- 1 }$. Together with Theorem \ref{theo:main}, applied
with $T = 1$, and with
$\mathcal W_2^2 (X_1, \mathbb X_1) \le 2 \mathcal W_2^2 (X_1, \widehat X_1)
+ 2 \mathcal W_2^2 ( \widehat X_1, \mathbb X_1 )$, this yields
\[
\mathcal W_2^2 ( X_1, \mathbb X_1 ) \le C \Big( N^{-p/2}
+ \sum_{r=1}^{m-1} \| c_N \otimes_r c_N \|_{m-r,2}^2 + n^{- 1 } \Big) ,
\]
which concludes the proof of Theorem \ref{theo:mainbis}.
\end{proof}

\begin{corollary}[Conditional propagation of chaos]\label{cor:chaos}
Grant the assumptions of Theorem \ref{theo:mainbis}  for some $q>4$ and suppose moreover that
$\E|X_0^1|^q<\infty$. Let
\[
\mu^N_1:=\frac1N\sum_{i=1}^N \delta_{X^{i}_1}
\qquad\text{and}\qquad
\mu_1:=\mathcal{L}\big(\mathbb{X}^1_1 \,\big|\, \sigma\{B_s,\, s\ge 0\}\big)
\]
be the empirical measure of the particle system at time $1$ and the conditional
law, given the common noise $B$, of one coordinate of the solution of (\ref{eq:SDE}).
Then it is possible to construct both systems on a common probability space  such that for every $p<p(q)$,
\[
\E\, W_2^2\big(\mu^N_1,\mu_1\big)
\le C\Big( N^{-p/2} + \sum_{r=1}^{m-1}\|c_N\otimes_r c_N\|^2_{m-r,2}
+ \frac1n + \frac1{\sqrt N}\Big),
\]
where $C$ depends on the model parameters and on $\E|X_0^1|^q$.
\end{corollary}

\begin{proof}
We briefly give the main ideas of the proof. 
The proofs of Theorems \ref{theo:main} and \ref{theo:mainbis}, based on
Proposition \ref{prop:41}, provide a coupling of $X_1$ with $(\mathbb{X}_1, B)$
on a common probability space such that $\E|X_1-\mathbb{X}_1|_N^2$ is upper bounded by
the right-hand side of Theorem \ref{theo:mainbis}. Write
$\bar\mu^N_1:=\frac1N\sum_{i=1}^N\delta_{\mathbb{X}^i_1}$. Since
$W_2^2\big(\frac1N\sum_i\delta_{x^i},\frac1N\sum_i\delta_{y^i}\big)\le |x-y|_N^2$
for all $x, y \in \R^N$ (see Remark \ref{rem:empirical}),
\[
\E\, W_2^2(\mu^N_1,\mu_1)
\le 2\,\E|X_1-\mathbb{X}_1|_N^2 + 2\,\E\, W_2^2(\bar\mu^N_1,\mu_1).
\]
Conditionally on $B$, the variables $\mathbb{X}^1_1,\dots,\mathbb{X}^N_1$ are i.i.d.
with common law $\mu_1$. Applying Theorem 1 of \cite{FG15} conditionally on $B$
(with $d=1$ and $q>4$) gives
\[
\E\big[W_2^2(\bar\mu^N_1,\mu_1)\,\big|\,B\big]
\le C\,\Big(\E\big[|\mathbb{X}^1_1|^q\,\big|\,B\big]\Big)^{2/q} N^{-1/2},
\]
and from this the assertion follows, taking expectation. 
\end{proof}

\begin{remark}
Corollary \ref{cor:chaos} states the conditional propagation of chaos property in
its usual form: the empirical measure of the finite system is close to the
conditional law of the limit given the common noise $B$, which is a random
probability measure. Since $p/2<1/6$, the term $N^{-1/2}$ coming from the
empirical measure does not change the rate of convergence.
\end{remark}

\section{Probabilistic tools: invariance principle, law of large
numbers and CLT}\label{sec:3}

\subsection{The invariance principle}
Throughout this section we work on some fixed probability space $ ( \Omega, \mathcal A, P) $ on which all random variables are defined. 

\subsubsection{Distances}

In what follows we will use the following distances between two probability measures $\mu 
$ and $\nu $ on ${\mathbb{R}}.$ For $k\in \{0,1,...\},$ we define 
\begin{equation*}
d_{k}(\mu ,\nu )=\sup \left\{\left\vert \int h d\mu -\int h d\nu
\right\vert :\sum_{i=0}^{k}\left\Vert h^{(i)}\right\Vert _{\infty
}\leq 1\right\}.
\end{equation*}

For $k=0,$ $ d_0 $ equals twice the total variation distance. For $k=1, $ $d_1$ is the
Fortet-Mourier distance. We will use $d_{3},$ which naturally appears when
using Lindeberg's method in connection with a Taylor development of order
three.

Moreover for probability measures on $ \R$ having a finite second moment we 
shall also use the quadratic Wasserstein distance $W_2$ given by 
\begin{equation*}
W_{2}^{2}(\mu ,\nu ):=\inf_{\pi \in \Pi _{\mu ,\nu }}\int_{{\mathbb{R}}%
^{2}}\left\vert x-y\right\vert ^{2}\pi (dx,dy),
\end{equation*}%
where $\Pi _{\mu ,\nu }$ is the set of probability measures on ${\mathbb{R}}^{2}$
having marginals $\mu $ and $\nu .$ By a slight abuse of notation, for two
random variables $X$ and $Y,$ we write for short $W_{2}(X,Y)$ instead of $ W_2 ( \mathcal L ( X), \mathcal L ( Y) ) $
for  the Wasserstein
distance of their laws. 

Finally, we recall that, when dealing with the $ N-$particle systems,
we work with vectors in ${\mathbb{R}}^{N}$ on which we consider the norm 
$$
\left\vert x\right\vert _{N}^{2}=\frac{1}{N}\sum_{i=1}^{N}\left\vert
x^{i}\right\vert ^{2},
$$
which induces the associated quadratic Wasserstein distance for probability measures on $ \R^N, $ denoted by $\mathcal{W}_{2}.$

\begin{remark}\label{rem:empirical}
Let us shortly comment on the reason why we have to work with $ \mathcal{W}_2.$
Notice that if $I$ is uniformly distributed on $\{1,\dots,N\}$ and independent of everything else, then $\mathcal W_2^2(X,Y)=\inf \E|X^I-Y^I|^2$, where the infimum is taken over all couplings of the random vectors $X$ and $Y$. So $\mathcal W_2$ measures the error for a typical particle, chosen at random. Since our system is not exchangeable, this is weaker than a bound for each fixed particle. However, since $W_2^2\big(\frac1N\sum_i\delta_{x^i},\frac1N\sum_i\delta_{y^i}\big)\le|x-y|_N^2$ for all $x,y\in\R^N$, it implies a bound for the empirical measures: for any coupling of $X$ and $Y$, $\E\,W_2^2\big(\frac1N\sum_i\delta_{X^i},\frac1N\sum_i\delta_{Y^i}\big)\le\E|X-Y|_N^2$, and the infimum of the right-hand side over all couplings is $\mathcal W_2^2(X,Y)$.
\end{remark}

As we have explained in the introduction, at some point we need to bound the
Wasserstein distance $W_{2}$ in terms of $d_{3}$ (see (\ref{W4}), (\ref{CLT13}) together with the conditions of Prop. \ref{prop:41} below). In order to do so, we rely
on Lemma B2 of \cite{ABCK}. For the convenience of the reader we recall its
statement here.

\begin{lemma}[Corollary of Lemma B2, \cite{ABCK}]\label{lem:21}
Let $\mu $ and $\nu $ be two probability measures on ${\mathbb{R}}.$ Then
for every $q>2$ and every 
\begin{equation*}
p<p(q):=\frac{2(q-2)}{6(q-2)+6}, 
\end{equation*}%
(recall (\ref{pq})),
we have 
\begin{equation}
W_{2}^{2}(\mu ,\nu )\leq C_{p}(1+\int \left\vert x\right\vert ^{q}\mu
(dx)+\int \left\vert x\right\vert ^{q}\nu (dx))^{1-p}d_{3}^{p}(\mu ,\nu ) .
\label{W1}
\end{equation}%
In particular, if $\int \left\vert x\right\vert ^{q}\mu (dx)+\int \left\vert
x\right\vert ^{q}\nu (dx)<\infty $ for every $q>2,$ then we have for every $%
p<\frac{1}{3},$ 
\begin{equation}
W_{2}^{2}(\mu ,\nu )\leq C_{p}d_{3}^{p}(\mu ,\nu ),  \label{W2}
\end{equation}%
for a constant $C_{p}$ depending only on $p$ and on the moments of order $ q $ of $ \mu $ and $ \nu, $ where $ p <  p (q) .$ 
\end{lemma}

\begin{proof}
The statement is a corollary of Lemma B2 of \cite{ABCK} which we recall here. For a set
of functions $\mathcal{A}$ we define 
\begin{equation*}
d_{\mathcal{A}}(\mu ,\nu )=\sup_{f\in \mathcal{A}}\left\vert \int fd\mu
-\int fd\nu \right\vert .
\end{equation*}%
Moreover, we choose a smooth probability density $\phi $ having support included
in the unit ball. Then the statement of Lemma B2 of \cite{ABCK} is as follows: 

If there exists $\kappa >0$ such that for every $%
x\in {\mathbb{R}}$ and for any $ r \in (0, 1 ] $ 
\begin{equation}\label{eq:kappa}
\left\vert \int \frac{1}{r}\phi (\frac{x-y}{r})d\mu (y)-\int \frac{1}{r}\phi
(\frac{x-y}{r})d\nu (y)\right\vert \leq \frac{C}{r^{\kappa }}d_{\mathcal{A}%
}(\mu ,\nu ),
\end{equation}%
then, putting 
\begin{equation}
p_{\kappa }(q)=\frac{2(q-2)}{(\kappa +2)(q-2)+6},  \label{general4}
\end{equation}%
we have that for every $p<p_{\kappa }(q)$ and for a convenient constant,   
$$%
W_{2}^{2}(\mu ,\nu )\leq C_{p}(1+\int \left\vert x\right\vert ^{q}\mu
(dx)+\int \left\vert x\right\vert ^{q}\nu (dx))^{1-p}d_{\mathcal{A}}^{p}(\mu ,\nu ).$$ 

We apply this result with  $\mathcal{A=}\{\varphi \in C_{b}^{3}(\R);\sum_{i=0}^{3}\left\Vert
\varphi ^{(i)}\right\Vert _{\infty }\leq 1\},$ so that $d_{\mathcal{A}%
}=d_{3}.$ It is immediate to check that inequality (\ref{eq:kappa}) holds with $\kappa =4.$ This
proves (\ref{W1}). It is immediate to obtain (\ref{W2}): for any $ p < \frac13, $ we choose $q$ such that $ p < p(q).$ 
\end{proof}

\begin{remark}
A similar estimate of $W_{2}$ in terms of the Kolmogorov distance has
recently been given by Pag\`es and Panloup \cite{PP26}. But, as we have
already mentioned, for our purpose, it is preferable to work with $d_{3}$ since 
this gives slightly more accurate bounds for the error.
\end{remark}

\subsection*{Multi-linear polynomials}
Recall that we work with a fixed number of particles $N,$ that $m  $ is the order of the interaction
and that $\Lambda _{N,m}$ denotes the set of multi-indices $\alpha
=(\alpha _{1},...,\alpha _{m})\in \{1,...,N\}^{m}.$ For any vector $%
x=(x^{1},...,x^{N})\in {\mathbb{R}}^{N},$ we denote $x^{\alpha
}=x^{\alpha _{1}}\times ...\times x^{\alpha _{m}}.$

We consider a collection of non-negative symmetric coefficients $c_{N}(\alpha
),\alpha \in \Lambda _{N,m},$ which are null on the diagonal, such that $ \|c_N\|_2 < \infty . $ We do not impose any 
further
normalization condition within this section.

Following \cite{NPR}, we consider \textit{multi-linear polynomials}
given by 
\begin{equation}
S_{N}(c_{N},x)=\sum_{\alpha \in \Lambda _{N,m}}c_{N}(\alpha )x^{\alpha
},\quad x\in {\mathbb{R}}^{N}.  \label{eq:polynom}
\end{equation}

We first state the following invariance principle for such polynomials which
allows us to switch from a random input $X $ to another input $Y$ and to
compare $S_N ( c_N, X) $ and $S_N ( c_N, Y) .$ This invariance principle is
due to Mossel et al. \cite{MDO}, see also \cite{NPR}. We use it in the version proved in \cite{BC}, Theorem 2.2.

\begin{proposition}
(\textbf{Invariance principle) }Let $X=(X^{1},...,X^{N})$ and $%
Y=(Y^{1},...,Y^{N})$ be two $N-$dimensional random vectors with independent
components such that $\E(X^{i})=\E(Y^{i})=0$ and $%
\E((X^{i})^{2})=\E((Y^{i})^2)=1$ for all $1 \le i \le N. $ Suppose
furthermore that they have finite moments of order three, that is, $%
\E\left\vert X^{i}\right\vert ^{3}<\infty $ and $\E\left\vert Y^{i}\right\vert
^{3}<\infty ,$ for all $1 \le i \le N.$ Finally, we suppose only that $ \|c_N\|_2^2 < \infty$ (no renormalization condition is imposed).
Then there exists a constant $K$
depending only on $\E\left\vert X^{i}\right\vert ^{3}$ and $\E\left\vert
Y^{i}\right\vert ^{3},i=1,...,N,$ and on $m,$  such that 
\begin{equation}
d_{3}(S_{N}(c_{N},X),S_{N}(c_{N},Y))\leq K \|c_N\|_2^2 \delta _{N}(c_{N}).  \label{ML1}
\end{equation}
\end{proposition}

\begin{remark}
In Mossel et al. \cite{MDO}, the authors provide an estimate in the
Kolmogorov distance. We prefer to work with the $d_{3}-$distance as in \cite%
{NPR} and in \cite{BC} which provides more accurate bounds on the speed of
convergence.
\end{remark}

We now want to combine the invariance principle (\ref{ML1}) and the estimate (\ref{W2}). To be able to do so, we need estimates for the moments of $S_N ( c_N, X) $ and $ S_N ( c_N, Y) $
which are a consequence of Proposition 3.11, Proposition
3.12 and Proposition 3.16 in \cite{MDO}, see also Lemma 4.2 in \cite{NPR}.

\begin{lemma}
(\textbf{Moments) }Let $Z^i, i \geq 1 , $ be independent centred random variables
with $\E\left\vert Z^{i}\right\vert ^{2}=1 $ for all $ i \geq 1.$ Suppose that for some $q>2$ one
has 
\begin{equation}
\gamma_Z (q):= \sup_{ i \geq 1 } \E\left\vert Z^{i}\right\vert ^{q}<\infty .  \label{gama}
\end{equation}%
Then 
\begin{equation}
\E\left\vert S_{N}(c_{N},Z)\right\vert ^{q}\leq
C_{q,m}\gamma_Z (q)^{m} (m! \|c_N\|_2^2)^{q/2} <\infty .  \label{ML4}
\end{equation}%
In particular, if $ \sup_N \|c_N\|_2 < \infty, $ then also $ \sup_N \E\left\vert S_{N}(c_{N},Z)\right\vert ^{q} < \infty .$ 
\end{lemma}

Recall that $p(q) $ was defined in (\ref{pq}). Together with (\ref{W2}) we therefore deduce from (\ref{ML1}) the following control.

\begin{proposition}
(\textbf{Invariance principle revisited) }
Let $q>2$ and suppose that $%
X=(X^{1},...,X^{N})$ and $Y=(Y^{1},...,Y^{N})$ are two $N-$dimensional
random vectors with independent components such that for all $ 1 \le i \le N, $  $\E(X^{i})=\E(Y^{i})=0,$ $%
\E((X^{i})^{2})=\E((Y^{i})^2)=1$ and $\sup_{ 1 \le i \le N} (\E|X^{i}|^{q}+\E|Y^{i}|^{q})=:\gamma_{q}<\infty .$ Then for every $p<p(q),$ 
\begin{equation}
W_{2}^{2}(S_{N}(c_{N},X),S_{N}(c_{N},Y))\leq C_{p,q} \|c_N\|_2^{ 2p}\delta _{N}^{p}(c_{N}),
\label{W4}
\end{equation}%
where $C_{p,q}$ depends on $p$ and on $\gamma_{q}. $ If $\gamma_{q}<\infty $ for every $%
q\in \N,$ then for every $p<\frac{1}{3},$ 
\begin{equation}
W_{2}^{2}(S_{N}(c_{N},X),S_{N}(c_{N},Y))\leq C_{p} \|c_N\|_2^{ 2p} \delta _{N}^{p}(c_{N}).
\label{W4'}
\end{equation}
\end{proposition}

\bigskip

\subsection{Conditional law of large numbers}
We now state a result that plays the role of the law of large numbers in our
framework. Let $F=(F^{1},\ldots ,F^{N})$ be an $N-$dimensional random
vector. We define 
\begin{equation}
\sigma _{N}^{2}(c_{N},F)=m!\sum_{\alpha \in \Lambda _{N,m}}c_{N}^{2}(\alpha
)F^{\alpha }.  \label{eq:cNx}
\end{equation}
Throughout this section we work again under our renormalization condition $ m! \|c_N\|_2^2 = 1.$ 

\begin{proposition}\label{prop:LLN}
\textbf{(Conditional law of large numbers)}\newline
Let $\mathcal{G}$\ be some $\sigma -$algebra and let $\E_{\mathcal{G}}$ be
the conditional expectation with respect to $\mathcal{G}.$\ Let $%
F=(F^{1},...,F^{N})$  be some $N-$dimensional random vector having square
integrable components which are $i.i.d.$ conditionally with respect to $\mathcal{G}.$\ 
Suppose moreover that%
\begin{equation*}
\mathcal{M}_{2,m}(F):=\E\left[ (\E_{\mathcal{G}}\left\vert F^{1}\right\vert
^{2})^{m}\right] <\infty .
\end{equation*}%
We write $\E_{\mathcal{G}}F=(\E_{\mathcal{G}}F^{1},...,\E_{\mathcal{G}}F^{N})$ for
the associated vector of conditional expectations. Then 
\begin{equation}
\E\left[ (\sigma _{N}^{2}(c_{N},F)-\sigma _{N}^{2}(c_{N},\E_{\mathcal{G}%
}F))^{2}\right] \leq m^2  (m!)^{4}\mathcal{M}_{2,m}(F)\widehat{\delta }%
_{N}^{2}(c_{N}).  \label{ML2}
\end{equation}%
In particular, since the components are identically distributed, one has 
\begin{equation*}
\sigma _{N}^{2}(c_{N},\E_{\mathcal{G}}F)=m!(\E_{\mathcal{G}}F^{1})^{m}\sum_{%
\alpha \in \Lambda _{N,m}}c_{N}^{2}(\alpha )=(\E_{\mathcal{G}}F^{1})^{m},
\end{equation*}%
so that 
\begin{equation}
\E\left[ (\sigma _{N}^{2}(c_{N},F)-(\E_{\mathcal{G}}F^{1})^{m})^{2}\right]
\leq m^2  (m!)^{4}\mathcal{M}_{2,m}(F)\widehat{\delta }_{N}^{2}(c_{N}).
\label{ML3}
\end{equation}
\end{proposition}

\begin{proof}
\textbf{Step 1} Let us write $\Lambda _{N,m}^{\uparrow }=\{\alpha \in
\Lambda _{N,m}:\alpha _{1}<\alpha _{2}<\ldots <\alpha _{m}\}.$ For $%
k=0,...,m,$ we denote%
\begin{equation*}
\sigma _{N,k}^{2}(c_{N},F)=(m!)^{2}\sum_{\alpha \in \Lambda _{N,m}^{\uparrow
}}c_{N}^{2}(\alpha )\Pi _{i=1}^{k}F^{\alpha _{i}}\times \Pi _{i=k+1}^{m}\E_{%
\mathcal{G}}F^{\alpha _{i}},
\end{equation*}%
so that $\sigma _{N,0}^{2}(c_{N},F)=\sigma _{N}^{2}(c_{N},\E_{\mathcal{G}}(F))$ and $\sigma
_{N,m}^{2}(c_{N},F)=\sigma _{N}^{2}(c_{N},F).$ Then 
\begin{equation*}
\E(\sigma _{N}^{2}(c_{N},F)-\sigma _{N}^{2}(c_{N},\E_{\mathcal{G}}F))^{2}\leq
m\sum_{k=1}^{m}\E(\sigma _{N,k}^{2}(c_{N},F)-\sigma
_{N,k-1}^{2}(c_{N},F))^{2}.
\end{equation*}%
We write%
\begin{eqnarray*}
&&\frac{1}{(m!)^{2}}\,\left( \sigma _{N,k}^{2}(c_{N},F)-\sigma
_{N,k-1}^{2}(c_{N},F)\right) \\
&=&\sum_{\alpha \in \Lambda _{N,m}^{\uparrow }}c_{N}^{2}(\alpha )\Pi
_{j=1}^{k-1}F^{\alpha _{j}}\times \left( F^{\alpha _{k}}-\E_{\mathcal{G}%
}F^{\alpha _{k}}\right) \times \Pi _{j=k+1}^{m}\E_{\mathcal{G}}F^{\alpha _{j}}
\\
&=&\sum_{i=1}^{N}\left( F^{i}-\E_{\mathcal{G}}F^{i}\right) \sum_{\alpha \in
\Lambda _{N,m}^{\uparrow },\alpha _{k}=i}c_{N}^{2}(\alpha )\Pi
_{j=1}^{k-1}F^{\alpha _{j}}\times \Pi _{j=k+1}^{m}\E_{\mathcal{G}}F^{\alpha
_{j}} \\
&=&\sum_{i=1}^{N}\left( F^{i}-\E_{\mathcal{G}}F^{i}\right) Q_{i}^{k},
\end{eqnarray*}%
where 
\begin{equation*}
Q_{i}^{k}=\sum_{\alpha \in \Lambda _{N,m}^{\uparrow },\alpha
_{k}=i}c_{N}^{2}(\alpha )\Pi _{j=1}^{k-1}F^{\alpha _{j}}\times \Pi
_{j=k+1}^{m}\E_{\mathcal{G}}F^{\alpha _{j}}.
\end{equation*}

\textbf{Step 2 }We estimate $\E_{\mathcal{G}}\left\vert Q_{i}^{k}\right\vert
^{2}.$ We write%
$$
\E_{\mathcal{G}}\left\vert Q_{i}^{k}\right\vert ^{2}\le \sum_{\alpha ,\beta \in
\Lambda _{N,m}^{\uparrow },\alpha _{k}=\beta _{k}=i}c_{N}^{2}(\alpha
)c_{N}^{2}(\beta )
\E_{\mathcal{G}}(\Pi _{j=1}^{k-1}|F^{\alpha _{j}}| |F^{\beta
_{j}}|)\times \Pi _{j=k+1}^{m}\E_{\mathcal{G}}|F^{\alpha _{j}}|\E_{\mathcal{G}%
}|F^{\beta _{j}}|.
$$
Notice that $\alpha _{1}<...<\alpha _{k-1}$ and $\beta _{1}<...<\beta _{k-1}$
so there are at most $k-1$ indices which are repeated in the concatenation $%
(\alpha _{1},...,\alpha _{k-1},\beta _{1},...,\beta _{k-1}).$ Moreover, if they
are repeated, then this is at most one time. Using also that $\E_{\mathcal{G}%
}|F^{\alpha _{j}}|  \leq (\E_{\mathcal{G}} |
F^{\alpha _{j}}|^{2})^{1/2},$ we obtain $\E_{\mathcal{G}}(\Pi
_{j=1}^{k-1}|F^{\alpha _{j}}| |F^{\beta _{j}}|)\leq (\E_{\mathcal{G}}\left\vert
F^{1}\right\vert ^{2})^{k-1}$ and finally%
\begin{eqnarray*}
\E_{\mathcal{G}}\left\vert Q_{i}^{k}\right\vert ^{2} &\leq &(\E_{\mathcal{G}%
}\left\vert F^{1}\right\vert ^{2})^{m-1}\sum_{\alpha ,\beta \in \Lambda
_{N,m}^{\uparrow },\alpha _{k}=\beta _{k}=i}c_{N}^{2}(\alpha
)c_{N}^{2}(\beta ) \\
&=&(\E_{\mathcal{G}}\left\vert F^{1}\right\vert ^{2})^{m-1}(\sum_{\alpha \in
\Lambda _{N,m}^{\uparrow },\alpha _{k}=i}c_{N}^{2}(\alpha ))^{2}.
\end{eqnarray*}

\textbf{Step 3 }Since the coordinates of $\alpha \in \Lambda
_{N,m}^{\uparrow }$ are strictly ordered, it follows that the variable $Q_{i}^{k}$ is $%
\sigma (\mathcal{G},F^{1},...,F^{i-1})-$measurable. Therefore, for $%
i_{1}<i_{2},$ conditionally to ${\mathcal{G}},$ $F^{i_{2}}$ is independent
of $\left( F^{i_{1}}-\E_{\mathcal{G}}F^{i_{1}}\right) Q^k_{i_{1}}Q^k_{i_{2}}.$ It
follows that%
\begin{equation*}
\E_{\mathcal{G}}\left( F^{i_{2}}-\E_{\mathcal{G}}F^{i_{2}} |
\sigma( F^{i_{1}}
,Q_{i_{1}}^{k},Q_{i_{2}}^{k},\mathcal{G}) \right) =\E_{\mathcal{G}}(
F^{i_{2}}-\E_{\mathcal{G}}F^{i_{2}}  )=0,
\end{equation*}%
such that 
$$
\E_{\mathcal{G}}\left[ \left( F^{i_{1}}-\E_{\mathcal{G}}F^{i_{1}}\right)
Q_{i_{1}}^{k}\times \left( F^{i_{2}}-\E_{\mathcal{G}}F^{i_{2}}\right)
Q_{i_{2}}^{k}\right] = 0.$$
Therefore, by orthogonality, 
\begin{multline*}
\E_{\mathcal{G}}(\sigma _{N,k}^{2}(c_{N},F)-\sigma _{N,k-1}^{2}(c_{N},F))^{2}
=(m!)^{4}\sum_{i=1}^{N}\E_{\mathcal{G}}\left( F^{i}-\E_{\mathcal{G}%
}F^{i}\right) ^{2}\times \E_{\mathcal{G}}|Q_{i}^{k}|^{2} \\
\leq (\E_{\mathcal{G}}\left\vert F^{1}\right\vert
^{2})^{m}(m!)^{4}\sum_{i=1}^{N}(\sum_{\alpha \in \Lambda _{N,m},\alpha
_{k}=i}c_{N}^{2}(\alpha ))^{2}.
\end{multline*}%
Taking expectation, we obtain%
\begin{multline*}
\E(\sigma _{N,k}^{2}(c_{N},F)-\sigma _{N,k-1}^{2}(c_{N},F))^{2} 
\leq \E((\E_{%
\mathcal{G}}\left\vert F^{1}\right\vert
^{2})^{m})(m!)^{4}\sum_{i=1}^{N}(\sum_{\alpha \in \Lambda _{N,m},\alpha
_{k}=i}c_{N}^{2}(\alpha ))^{2} \\
\le (m!)^{4}\mathcal{M}_{2,m}(F)\widehat{\delta }_{N}^{2}(c_{N}),
\end{multline*}
which allows to conclude. 
\end{proof}

\subsection{Coming back to the discrete time model}
We come back to the discrete time model of Section \ref{sec:frameworkgeneral}. We fix $ s < t $ and consider the associated
family of random fields $ v^i_{s, t} (x^i) $ defining the perturbations. Since we work for fixed time within this section, in what follows, 
we drop all time indices and write for short $ v^i (x^i ) $ instead of $ v^i_{s, t} ( x^i), $ etc. 
So we work with the family of random fields $v^{i}(x^i)=v^{i}(x^i,\omega ),$ $1 \le i \le N, $ $x^i \in \R ,$ satisfying Assumption
\ref{ass:perturbation} and Assumption \ref{ass:2.10first} for some $ q \geq 3.$ Again, throughout this section we impose the normalization conditions 
$$ m! \|c_N\|_2^2 = 1 \mbox{ and }  \| \bar c_N\|_2^2 = 1.$$

For any $\alpha =(\alpha _{1},...,\alpha _{m})\in \Lambda _{N,m}$ and $%
x=(x^{1},...,x^{N})\in {\mathbb{R}}^{N},$ we put%
\begin{equation*}
v^{\alpha }(x)=v^{\alpha _{1}}(x^{\alpha _{1}})\times ...\times v^{\alpha
_{m}}(x^{\alpha _{m}}).
\end{equation*}%
Recall (\ref{eq:polynom}). With this notation, we
associate to the vector $v(x)=(v^{1}(x^{1}),...,v^{N}(x^{N}))$ the random variable
\begin{equation}
S_{N}(c_{N},v(x))=\sum_{\alpha \in \Lambda _{N,m}}c_{N}(\alpha )v^{\alpha
}(x) \label{eq:S_Nx}
\end{equation}%
having variance
\begin{equation*}
\sigma _{N}^{2}(x)=\E (S_{N}^{2}(c_{N},x))=m! \sum_{\alpha }c_{N}^{2}(\alpha
)\prod_{i=1}^{m}f^{2}(x^{\alpha _{i}}) . 
\end{equation*}%
Notice that with the notation of  (\ref{eq:cNx}), $\sigma _{N}^{2}(x)=\sigma
_{N}^{2}(c_{N},F)$ with $F^{i}=f^2(x^{i}).$

First we give an application of the law of large Numbers. Recall that we work under the hypothesis that $ m! \|c_N\|_2^2 = 1.$

\begin{proposition}
Let $\mathcal{G} $ be a sub-$\sigma -$algebra of $ \mathcal A $ and let $X^1, \ldots, X^N $ be i.i.d., conditionally on $ \mathcal G.$ 
Then 
\begin{equation}
\E\left\vert \sigma _{N}(X)-(\E_{\mathcal{G}}f^{2}(X^{1}))^{m/2}\right\vert
^{2}\leq m^2 \times (m!)^{4}\underline{f}^{-2m}\overline{f}^{4m}\times \widehat{%
\delta }_{N}^{2}(c_{N}).  \label{CLT4}
\end{equation}
\end{proposition}

\begin{proof}
We have that $\E_{\mathcal{G}}(f^{2}(X^{1}))^{m}\geq \underline{f}^{2m}$ and $ \sigma_N^2 ( X) \geq  \underline{f}^{2m}$ as well. 
Therefore, using the identity $a-b=\frac{a^{2}-b^{2}}{a+b},$ 
\begin{multline*}
\E\left( \left\vert \sigma _{N}(X)-(\E_{\mathcal{G}}f^{2}(X^{1}))^{m/2}%
\right\vert ^{2}\right)  \\
\leq \underline{f}^{-2m}\E \left( \left\vert \sigma _{N}^{2}(X)-(\E_{%
\mathcal{G}}f^{2}(X^{1}))^{m}\right\vert ^{2}\right)  \\
=(m!)^2 \underline{f}^{-2m}\E\left\vert \sum_{\alpha }c_{N}^{2}(\alpha
) \left( \prod_{i=1}^{m}f^{2}(X^{\alpha _{i}})-(\E_{\mathcal{G}}f^{2}(X^{1}))^{m} \right) %
\right\vert ^{2}\\
\leq m^2 \times m!^{4}\underline{f}^{-2m}\overline{f}^{4m}%
\widehat{\delta }_{N}^{2}(c_{N}),
\end{multline*}%
the last inequality being a consequence of (\ref{ML3}), applied to  $%
F^{i}=f^{2}(X^{i}).$
\end{proof}

\subsection*{Gaussian representation}
For every fixed $x\in \mathbb{R}^N$ we define the coefficients 
\begin{equation*}
c_{N}(\alpha ,x)=c_{N}(\alpha )f(x^{\alpha _{1}})...f(x^{\alpha _{m}}).
\end{equation*}%
Moreover, we consider a sequence of independent standard normal random variables $%
G^{i},i=1,...,N,$ and we associate  the sum%
\begin{equation*}
S_{N}(c_{N}(x),G)=\sum_{\alpha \in \Lambda _{N,m}}c_{N}(\alpha ,x)G^{\alpha
}=\sum_{\alpha \in \Lambda _{N,m}}c_{N}(\alpha )f(x^{\alpha
_{1}})...f(x^{\alpha _{m}})G^{\alpha },
\end{equation*}%
with $G^{\alpha }=\prod_{i=1}^{m}G^{\alpha _{i}}.$ As an immediate
consequence of the invariance principle from the previous section we get

\begin{lemma} Suppose that Assumptions \ref{ass:perturbation} and \ref{ass:2.10first} hold for some $ q \geq 3$  for the random field $ v^i (x^i), 1 \le i \le N.$ Then for each fixed $x\in {\mathbb{R}}^{N},$%
\begin{equation}
d_{3}(S_{N}(c_{N}(x),G),S_{N}(c_{N},v(x)))\leq C\overline{f}^{3m}\delta
_{N}(c_{N}).  \label{CLT1}
\end{equation}%
Moreover for any $%
p<p(q)$ (recall (\ref{pq})), 
\begin{equation}
W^2_{2}(S_{N}(c_{N}(x),G),S_{N}(c_{N},v(x)))\leq C \overline{f}^{3mp} \delta
_{N}^{p}(c_{N}) .  \label{CLT1'}
\end{equation}%
The constant $ C$ appearing in (\ref{CLT1}) and (\ref{CLT1'}) depends in both cases on $\gamma(3), $ on $ \gamma ( q), \underline f$ and $ \overline f.$ 
In particular, if $\gamma (q)<\infty $ for every $q\in \N ,$ then (\ref{CLT1'}) holds for every $p<%
\frac{1}{3}.$ 
\end{lemma}

\begin{proof}
Let%
\begin{equation*}
Z(x)=(Z^{1}(x^{1}),...,Z^{N}(x^{N}))\quad \mbox{with} \quad Z^{i}(x^{i})=\frac{%
v^{i}(x^{i})}{f(x^{i})}.
\end{equation*}%
Then we may rewrite 
\begin{equation*}
S_{N}(c_{N},v(x))=\sum_{\alpha \in \Lambda _{N,m}}c_{N}(\alpha ,x)Z^{\alpha
}(x)=S_{N}(c_{N}(x),Z(x)).
\end{equation*}

Clearly, for every fixed $x=(x^{1},\ldots ,x^{N}),$ $Z^{i}(x^{i}),i=1,...,N,$
are independent, centred random variables with variance equal to one and with third moments upper bounded by $ \gamma (3) / \underline f^3.$  We apply the invariance principle (\ref{ML1}) of the previous section to $ Z^i ( x^i ) $ and to the coefficients $ c_N ( x).$ 
Then we get 
\begin{multline*}
d_{3}(S_{N}(c_{N},v(x)),S_{N}(c_{N}(x),G))=d_{3}(S_{N}(c_{N}(x),Z(x)),S_{N}(c_{N}(x),G))\\
\leq K \|c_N ( x) \|_2^2 \delta _{N} (c_{N}(x)).
\end{multline*}
Since for every $x\in {\mathbb{R}}^{N}$ one has 
\begin{equation*}
\delta _{N}^{2}(c_{N}(x))\leq \overline{f}^{2m}\delta _{N}^{2}(c_{N})
\end{equation*}%
and 
$$ \|c_N ( x) \|_2^2 \le \overline f^{2m} \|c_N\|_2^2 = \overline f^{2m} / m! , $$
this implies (\ref{CLT1}). The second assertion is an immediate consequence of %
(\ref{W2}). 
\end{proof}

Once we may deal with $S_{N}(c_{N}(x),G),$ we are in the framework of the Wiener chaos. Indeed, since $c_N$ vanishes on the diagonal, $S_{N}(c_{N}(x),G)$
represents a multiple stochastic integral of order $m,$ for a piecewise
constant function. To be more precise, we construct the following function.
For any $t\in ] 0,1[,$ we define $q_{N}(t)=i+1,$ if $\frac{i}{N}\leq t<%
\frac{i+1}{N}.$ Moreover, for every fixed $x=(x^{1},...,x^{N}),$ we define%
\begin{equation}
\phi
_{N}(x^{1},...,x^{N};t_{1},...,t_{m})=N^{m/2}c_{N}(q_{N}(t_{1}),..,q_{N}(t_{m}))f(x^{q_{N}(t_{1})})...f(x^{q_{N}(t_{m})}).
\label{CLT9'}
\end{equation}%
In what follows we represent $ G^1, \ldots, G^N $ by means of a standard Brownian motion $W$ via
$$G^{i}=\sqrt{N}(W(\frac{i}{N})-W(%
\frac{i- 1}{N})),i=1,...,N,$$   
such that we may rewrite, for each $%
x=(x^{1},...,x^{N}),$ 
\begin{eqnarray}
S_{N}(c_{N}(x),G) &=&\int_{0}^{1}...\int_{0}^{1}\phi
_{N}(x^{1},...,x^{N}; t_{1},...,t_{m})dW_{t_{1}} \ldots d W_{t_{m}}  \label{CLT9} \\
&=&I_{m}(\phi _{N}(x; \circ )) .  \notag
\end{eqnarray}%
Here we use the convention of \cite{NP}: for symmetric $ \phi \in L^2 ( [ 0, 1]^m),$ $ I_m ( \phi) $ denotes the multiple Wiener-It\^o integral over $ [0, 1]^m $ such that 
$ \E I_m(\phi)^2 = m! \| \phi\|^2 .$ 

\subsection{The CLT on the Wiener chaos}

Let $W=(W_u)_{u\in[0,1]}$ be an auxiliary standard
Brownian motion and, for a symmetric $\varphi\in L^2([0,1]^m)$, let $%
I_m(\varphi)$ be the multiple Wiener--It\^o integral of order $m$ with
respect to $W$, so that ${\mathbb{E}} I_m(\varphi)^2=m!\|\varphi%
\|_{L^2([0,1]^m)}^2$. For $1\leq r\leq m-1$ the contraction $%
\varphi\otimes_r\varphi\in L^2([0,1]^{2(m-r)})$ is 
\begin{equation*}
\varphi\otimes_r\varphi(s_1,\ldots,s_{m-r},u_1,\ldots,u_{m-r})=%
\int_{[0,1]^r}\varphi(s,v)\varphi(u,v)\,dv . 
\end{equation*}

\begin{theorem}[CLT, \protect\cite{NuPe,NP,LNP}]
\label{theo:CLT} Let $\varphi\in L^2([0,1]^m)$ be symmetric, $F=I_m(\varphi)$
and $v^2={\mathbb{E}} F^2>0$. Then 
\begin{equation}  \label{CLT10}
W_2^2(F,v\Delta)\leq\frac{C_m}{v^2}\sum_{r=1}^{m-1}\|\varphi\otimes_r\varphi%
\|^2_{L^2([0,1]^{2(m-r)})}.
\end{equation}
\end{theorem}

\begin{proof}
We quickly explain how the above result follows from the literature. A Stein kernel of a
centred random variable $F^{\prime }$ is a random variable $\tau _{F^{\prime
}}(F^{\prime })$ such that 
$${\mathbb{E}}(F^{\prime }g(F^{\prime }))={\mathbb{%
E}}(g^{\prime }(F^{\prime })\tau _{F^{\prime }}(F^{\prime }))$$ 
for all
smooth $g.$ Moreover, the Stein discrepancy is 
$$S^{2}(F^{\prime })={\mathbb{E}}%
\left\vert 1-\tau _{F^{\prime }}(F^{\prime })\right\vert^2 .$$ By \cite[%
Proposition 3.1]{LNP}, $W_{2}(F^{\prime },\Delta )\leq S(F^{\prime })$. We now take 
$F^{\prime }=F/v$. The duality formula of Malliavin calculus then gives 
$$\tau
_{F^{\prime }}(F^{\prime })={\mathbb{E}}(m^{-1}\left\Vert DF^{\prime
}\right\Vert ^{2}\mid F^{\prime }),$$ where $\Vert \cdot \Vert $ is the norm
of $L^{2}([0,1])$, and ${\mathbb{E}}(m^{-1}\left\Vert DF^{\prime
}\right\Vert ^{2})={\mathbb{E(}}F^{\prime })^{2}=1$. Hence, by Jensen's
inequality, 
\begin{equation*}
S^{2}(F^{\prime })\leq \mathrm{Var}\big(m^{-1}\left\Vert DF^{\prime
}\right\Vert ^{2}\big)=v^{-4}\,\mathrm{Var}\big(m^{-1}\Vert DF\Vert ^{2}\big)%
\leq C_{m}v^{-4}\sum_{r=1}^{m-1}\Vert \varphi \otimes _{r}\varphi \Vert ^{2}.
\end{equation*}%
The last inequality is a consequence of the product formula for multiple
integrals, see \cite{NP} or \cite[Chapter 5]{NPbook}. Finally $%
W_{2}^{2}(F,v\Delta )=v^{2}W_{2}^{2}(F^{\prime },\Delta )$.
\end{proof}

For a fixed $x=(x^{1},...,x^{N}),$ we apply the above result to
\begin{equation*}
\varphi (t_{1},...,t_{m})=\phi _{N}(x^{1},...,x^{N};t_{1},...,t_{m})\times
\left\Vert I_{m}(\phi _{N}(x,\circ ))\right\Vert _{2}^{-1},
\end{equation*}%
where $\phi _{N}$ was defined in (\ref{CLT9'}). Since the $ c_N ( \alpha ) $ are supposed to be non-negative, 
it is straightforward to check that for every $x,$
we have
\begin{equation*}
\left\Vert \phi _{N}(x,\circ )\otimes _{r}\phi _{N}(x,\circ )\right\Vert
_{m-r}^{2}\leq \overline{f}^{4m}\left\Vert c_{N}\otimes _{r}c_{N}\right\Vert
_{m-r}^{2}.
\end{equation*}%
This together with the previous estimate (\ref{CLT10}) yields 
\begin{eqnarray*}
&&W_{2}^{2}(I_{m}(\phi _{N}(x,\circ )),\left\Vert I_{m}(\phi _{N}(x,\circ
))\right\Vert _{2}\Delta )   \\
&&= \left\Vert I_{m}(\phi _{N}(x,\circ ))\right\Vert _{2}^{2}W_{2}^{2}(\frac{%
I_{m}(\phi _{N}(x,\circ ))}{\left\Vert I_{m}(\phi _{N}(x,\circ ))\right\Vert
_{2}},\Delta )\\
 &&=\left\Vert I_{m}(\phi _{N}(x,\circ ))\right\Vert
_{2}^{2}W_{2}^{2}(I_{m}(\varphi ),\Delta ) \\
&&\leq C_{m}\left\Vert I_{m}(\phi _{N}(x,\circ ))\right\Vert
_{2}^{2}\sum_{r=1}^{m-1}\left\Vert \varphi \otimes _{r}\varphi \right\Vert
_{m-r}^{2} \\
&&=C_{m} \left\Vert I_{m}(\phi _{N}(x,\circ ))\right\Vert
_{2}^{- 2}\sum_{r=1}^{m-1}\left\Vert \phi _{N}(x,\circ )\otimes _{r}\phi
_{N}(x,\circ )\right\Vert _{m-r}^{2} \\
&&\leq C_{m}\overline{f}^{4m}/\underline f^{2m}\sum_{r=1}^{m-1}\left\Vert c_{N}\otimes
_{r}c_{N}\right\Vert _{m-r}^{2},
\end{eqnarray*}%
where we have used that $\left\Vert I_{m}(\phi _{N}(x,\circ ))\right\Vert
_{2}^{2}\geq \underline f^{2m}.$ 

Recall that $I_{m}(\phi _{N}(x,\circ ))=S_{N}(c_{N}(x),G)$ and $\sigma
_{N}^{2}(x)=\left\Vert S_{N}(c_{N}(x),G)\right\Vert _{2}^{2}=\left\Vert
I_{m}(\phi _{N}(x,\circ ))\right\Vert _{2}^{2}.$ So all in all we get%
\begin{equation}\label{eq:barepsilon}
W_{2}^{2}(S_{N}(c_{N}(x),G),\sigma _{N}(x)\Delta )\leq  \varepsilon_N( c_N) :=  C_{m}\overline{f}^{4m}/\underline f^{2m}%
\sum_{r=1}^{m-1}\left\Vert c_{N}\otimes _{r}c_{N}\right\Vert _{m-r, 2 }^{2}.
\end{equation}%
Combining this with the invariance principle (\ref{W4}) we conclude with the following estimate.

\begin{proposition}
Suppose that Assumptions \ref{ass:perturbation} and  \ref{ass:2.10first} hold for some $ q \geq 3.$ Then for every $p<p(q),$ 
\begin{equation}
W_{2}^{2}(S_{N}(c_{N},v(x)),\sigma _{N} (x)\Delta )\leq 2 C_{p, q }\delta
_{N}^{p}(c_{N})+2 C_{m}\overline{f}^{4m}/\underline f^{2m}\sum_{r=1}^{m-1, 2 }\left\Vert
c_{N}\otimes _{r}c_{N}\right\Vert _{m-r}^{2}.  \label{CLT13}
\end{equation}%
If Assumption \ref{ass:2.10first} holds for every $q,$ then the above estimate holds for
every $p<\frac{1}{3}.$
\end{proposition}

\section{Proof of Theorem \ref{theo:main}}\label{sec:4}
We are now able to give the proof of Theorem \ref{theo:main}. The first main ingredient is a control of the Wasserstein distance 
$\mathcal{W}_{2}^{2}(X_{t_{k}},\overline{X}_{t_{k}})$ of the original dynamic and its first approximation. 
This control is obtained using the results of the previous section together with a general statement of \cite{EvaVlad} that we cite for the 
convenience of the reader.  

\begin{proposition}[Prop. 4.2 of \cite{EvaVlad}]\label{prop:41}
Suppose that the following Lipschitz properties hold:
\begin{equation}\label{eq:LipT}
\left\vert \mathbb{T}(x)-\mathbb{T}(y)\right\vert \leq  L_{\mathbb{T}}\left\vert x-y\right\vert _{N}  
\end{equation}
and
\begin{equation}
\E\left\vert \overline{S}_{s,t}^{N}(c_N, x )-\overline{S}_{s,t}^{N}(%
c_N, y )\right\vert ^{2}\leq \overline{L}(t-s) |x - y
|_N^2    \label{Aiii}
\end{equation}
for all $ x, y \in \R^N.$ Suppose furthermore that for all $x \in {\mathbb{%
R}}^{N}, $ we have the coupling property
\begin{equation}
W_{2}^{2}(S_{s,t}^{N}(c_N, x ),\overline{S}_{s,t}^{N}(c_N, x ))\leq (t-s)%
\overline{\varepsilon }_{N} ,  \label{Aiv}
\end{equation}%
where $\lim_{ N \to \infty } \bar \varepsilon_N = 0.$

Then there exists a constant $ Q$ such that
\begin{equation}
\mathcal{W}_{2}^{2}(X_{t_{k}},\overline{X}_{t_{k}})\leq
e^{Q t_{k}}(\mathcal{W}_{2}^{2}(X_{0},\overline{X}%
_{{0}})+2 t_{k}\overline{\varepsilon }_{N}) , \label{clt4}
\end{equation}%
and $Q$ depends only on $ L_{\mathbb{T}} , $ $\overline{L},$ $%
C_{stb} $ and the Lipschitz constant of $ f^2. $ 
\end{proposition}

\begin{remark}
The only difference with Proposition 4.2 of \cite{EvaVlad}  is that $\frac1N \sum_{i=1}^N \varphi ( x^i) $ is replaced by $\mathbb{T}(x).$
But the only property of $x \mapsto\frac1N \sum_{i=1}^N \varphi ( x^i)$ which was used in the proof of 
Prop. 4.2 in \cite{EvaVlad} is the Lipschitz property (\ref{eq:LipT}). 
Apart from this, the proof is identical, and so we skip it.
\end{remark}

\subsection{Checking the conditions of Proposition \ref{prop:41}} 
In this section we show that under the general conditions of Theorem \ref{theo:main}, the conditions used in Proposition \ref{prop:41}
are all satisfied. 

Recall that we suppose that for some constant $C > 0, $ 
\begin{equation}\label{eq:recall}
\delta _{N}({c}_{N})+ \delta _{N}(\overline{c}_{N})\leq \frac{C}{\sqrt{N}}. 
\end{equation}
We first show that this condition implies the Lipschitz property of $ \mathbb{T}.$

\begin{proposition}\label{prop:LipT}
Suppose that (\ref{McKV6}) holds. Then
\begin{equation}
\left\vert \mathbb{T}(x)-\mathbb{T} (y)\right\vert \leq C \left\Vert \varphi \right\Vert _{\infty }^{\overline{m}-1}L_{\varphi
}\left\vert x-y\right\vert _{N}.  \label{McKV4}
\end{equation}
In particular, (\ref{eq:LipT}) holds with $ L_{\mathbb{T}} = C \left\Vert \varphi \right\Vert _{\infty }^{\overline{m}-1}L_{\varphi
} . $
\end{proposition}

\begin{proof}
We write 
\begin{eqnarray*}
\left\vert \mathbb{T} (x)-\mathbb{T} (y)\right\vert  &\leq
&\sum_{\alpha \in \Lambda _{N,\overline{m}}}\overline{c}_{N}^{2}(\alpha
)\left\vert \prod_{j=1}^{\overline{m}}\varphi (x^{\alpha _{j}})-\prod_{j=1}^{%
\overline{m}}\varphi (y^{\alpha _{j}})\right\vert  \\
&\leq &\left\Vert \varphi \right\Vert _{\infty }^{\overline{m}-1}L_{\varphi
}\sum_{j=1}^{\overline{m}}\sum_{\alpha \in \Lambda _{N,\overline{m}}}%
\overline{c}_{N}^{2}(\alpha )\left\vert x^{\alpha _{j}}-y^{\alpha
_{j}}\right\vert .
\end{eqnarray*}%
Since $\delta _{N}(\overline{c}_{N})\leq C N^{-1/2},$ we get%
\begin{equation*}
\sum_{k=1}^{N}(\sum_{ \alpha \in \Lambda _{N,\bar m}^{k}}\overline{c}_{N}^{2}(\alpha ))^{2}\leq
N\delta _{N}^{4}(\overline{c}_{N})\leq C N^{-1}
\end{equation*}%
and consequently, 
\begin{eqnarray}
\sum_{j=1}^{\bar m}\sum_{\alpha \in \Lambda _{N,\overline{m}}}\overline{c}_{N}^{2}(\alpha
)\left\vert x^{\alpha _{j}}-y^{\alpha _{j}}\right\vert 
&=&\sum_{k=1}^{N} \left( \sum_{j=1}^{\bar m} \sum_{\alpha : \alpha_j = k } \overline{c}_{N}^{2}(\alpha
)\right)\left\vert x^{k}-y^{k}\right\vert   \nonumber  \\
&=& \sum_{k=1}^{N}\left(\sum_{\alpha \in \Lambda _{N,\bar m}^{k} }\overline{c}_{N}^{2}(\alpha
)\right)\left\vert x^{k}-y^{k}\right\vert  \nonumber \\
&\leq &\left(\sum_{k=1}^{N}\left(\sum_{ \alpha \in \Lambda _{N,\bar m}^{k}}\overline{c}_{N}^{2}(\alpha
)\right)^2\right)^{1/2}(\sum_{k=1}^{N}\left\vert x^{k}-y^{k}\right\vert ^{2})^{1/2}  \notag
\\
&\leq &(C N^{-1}\sum_{k=1}^{N}\left\vert x^{k}-y^{k}\right\vert
^{2})^{1/2}=C  \left\vert x-y\right\vert _{N}  \label{McKV4'} 
\end{eqnarray}%
(recall that constants are allowed to change from one occurrence to another). 
\end{proof}

We next check the Lipschitz property (\ref{Aiii}). 

\begin{proposition}\label{prop:LipS}
Grant Assumption \ref{ass:perturbation} and suppose that (\ref{eq:recall}) holds. Then, using the same Gaussian increment to represent  $ \overline{S}_{s,t}^{N}(c_{N},x) $ and 
$ \overline{S}_{s,t}^{N}(c_{N},y),$ we have 
\begin{equation}
\E\left( \left\vert \overline{S}_{s,t}^{N}(c_{N},x)-\overline{S}%
_{s,t}^{N}(c_{N},y)\right\vert ^{2}\right) \leq C (t-s)(m!)^2 \underline{f}^{-2m}%
\overline{f}^{4m-4}L_{f^2}^{2}\left\vert x-y\right\vert _{N}^{2}.
\label{McKV9'}
\end{equation}%
In particular, the Lipschitz property (\ref{Aiii}) holds for
fixed $s<t$\ with $$\bar{L}=C (m!)^2 \underline{f}^{-2m}\overline{f}^{4m-4}L_{f^2}^{2}.$$
\end{proposition}

\begin{proof}
We use the inequality 
\begin{equation*}
\left\vert \sigma (x)-\sigma (y)\right\vert =\left\vert \sigma
^{2}(x)-\sigma ^{2}(y)\right\vert \times (\sigma (x)+\sigma (y))^{-1}\leq 
\underline{f}^{-m}\left\vert \sigma ^{2}(x)-\sigma ^{2}(y)\right\vert 
\end{equation*}%
and obtain 
$$
\E\left\vert \overline{S}_{s,t}^{N}(c_{N},x)-\overline{S}%
_{s,t}^{N}(c_{N},y)\right\vert ^{2} 
=(t-s)\left\vert \sigma (x)-\sigma (y)\right\vert ^{2}\leq (t-s)\underline{%
f}^{-2m}\left\vert \sigma ^{2}(x)-\sigma ^{2}(y)\right\vert ^{2} .
$$
Moreover, 
\begin{eqnarray*}
\left\vert \sigma ^{2}(x)-\sigma ^{2}(y)\right\vert ^{2}
&\le  &(m!)^2 \left( \sum_{\alpha }c_{N}^{2}(\alpha
)\left\vert \prod_{i=1}^{m}f^{2}(x^{\alpha
_{i}})-\prod_{i=1}^{m}f^{2}(y^{\alpha _{i}})\right\vert \right) ^{2} \\
&\leq &(m!)^2 \overline{f}^{4m-4}L_{f^2}^{2}\left(
\sum_{i=1}^{m}\sum_{\alpha }c_{N}^{2}(\alpha )\left\vert x^{\alpha
_{i}}-y^{\alpha _{i}}\right\vert \right) ^{2} \\
&\leq &C (m!)^2 \overline{f}^{4m-4}L_{f^2}^{2}\left\vert
x-y\right\vert _{N}^{2},
\end{eqnarray*}%
where we have used (\ref{McKV4'}), applied to $c_N ,$  to obtain the last inequality.
\end{proof}

The core result that makes our method work is the following coupling estimate. 
It is a consequence of the arguments of Section \ref{sec:3} and follows from 
the invariance principle of Mossel et al. \cite{MDO}  and from a quantitative version of the Central Limit Theorem for multiple integrals, due to Nourdin
and Peccati  \cite{NP} and to Ledoux, Nourdin and 
Peccati \cite{LNP}.

\begin{proposition}\label{prop:coupling}
Grant Assumptions \ref{ass:perturbation} and  \ref{ass:2.10first} for some $q\geq 3.$ Then for every 
$$
p<p(q):=\frac{2(q-2)}{6(q-2)+6},  
$$
it holds that 
\begin{equation}
W_{2}^{2}(S_{s,t}^{N}(c_{N},x),\overline{S}_{s,t}^{N}(c_{N},x))\leq
C(t-s)(\delta _{N}^{p}(c_{N})+C_{m}\overline{f}^{4m}/\underline{f}^{2m}\sum_{r=1}^{m-1}\left%
\Vert c_{N}\otimes _{r}c_{N}\right\Vert _{m-r, 2 }^{2}),   \label{McKV10}
\end{equation}%
with $C$ depending on  $p, q,$ on $\gamma_T (q) ,$ on $\overline{f}$ and on $\underline{f%
}.$  
\end{proposition}

\begin{proof}
The proof is a direct consequence of (\ref{CLT13}), since $$ W_{2}^{2}(S_{s,t}^{N}(c_{N},x),\overline{S}_{s,t}^{N}(c_{N},x)) = 
(t-s) W_{2}^{2}(S_N (c_N, v (x) ), \sigma_N ( x) \Delta).
$$
\end{proof} 

As a consequence,  the coupling property (\ref{Aiv})
holds for fixed $s<t$ with 
\begin{equation}\label{eq:barepsilonN}
\overline{\varepsilon }_{N, p }= C \left( \delta _{N}^{p}(c_{N})+C_{m}\overline{f}^{4m}/\underline{f}^{2m}\sum_{r=1}^{m-1}\left%
\Vert c_{N}\otimes _{r}c_{N}\right\Vert _{m-r, 2 }^{2} \right) .
\end{equation}

Rephrasing Proposition 4.2 of \cite{EvaVlad} in our present context then gives

\begin{proposition}
Under the conditions of Theorem  \ref{theo:main} it holds that 
for every $p<p(q)$ and for every $ k $ with $t_k \le T,$  
\begin{equation}
\mathcal{W}_{2}^{2}(X_{t_{k}},\overline{X}_{t_{k}})\leq e^{Q T  }\left( 
\mathcal{W}_{2}^{2}(X_{0},\overline{X}_{0})+2 t_{k} \overline{\varepsilon }_{N, p } \right) ,  \label{CLT5}
\end{equation}%
with $\overline{\varepsilon }_{N, p }$ given in (\ref{eq:barepsilonN}) above.

If moreover $\max_{k : t_k \le T } \gamma_{t_{k-1}, t_k} (q) < \infty $ (recall (\ref{CLT2first})) for every $%
q\geq 3,$ then the above estimate holds for every $p<\frac{1}{3}.$
\end{proposition}

\subsection{Comparing the auxiliary chain to the limit dynamic}

The goal of this section is to compare the auxiliary chain to the limit chain. This
is mainly done by using the conditional law of large numbers that we have already 
established. Recall that we suppose that $ \delta_N ( c_N) + \delta_N ( \bar c_N) \le C/ \sqrt{N}.$ 

\begin{proposition}
\textbf{A }Let $X=(X^{1},\ldots ,X^{N})\in L^{2}_{t_k} (\R^N) $ such that $X^{1},\ldots ,X^{N}$ are 
i.i.d. conditionally on $\mathcal{F}_{t_{k}}^{\Delta }.$ Then for some generic constant $C_{\bar m} $ depending on $ \bar m, $ 
\begin{eqnarray}
\E\left\vert \E_{t_{k}}(h_{t_{k},t_{k+1}}(X)-\widehat{h}_{t_{k},t_{k+1}}(X))%
\right\vert _{N}^{2} &\leq &\frac{C_{\bar m} C_{stb}^{2}(t_{k+1}-t_{k})^{2}\left\Vert \varphi \right\Vert _{\infty
}^{2\overline{m}}}{N},  \label{McKV16} \\
\E\left\vert \widetilde{\E}_{t_{k}}(h_{t_{k},t_{k+1}}(X)-\widehat{h}%
_{t_{k},t_{k+1}}(X))\right\vert _{N}^{2} &\leq &\frac{C_{\bar m} C_{stb}^{2}(t_{k+1}-t_{k}) \left\Vert \varphi
\right\Vert _{\infty }^{2\overline{m}}}{N}.  \label{McKV17}
\end{eqnarray}%
\textbf{B }Recall that $\sigma _{N}^{2}(X)=m! \sum_{\alpha
}c_{N}^{2}(\alpha )\prod_{j=1}^{m}f^2(X^{\alpha _{j}}).$ Using the law of
large numbers (\ref{CLT4}), we obtain
\begin{equation}
\E\left[ \left| \sigma _{N}(X)-(\E_{t_{k}}^{\Delta }f^{2}(X^{1}))^{m/2}%
\right|^2 \right] \leq C_{m} \underline f^{ - 2m} \overline f ^{4m} N^{ - 1 } .  
\end{equation}
\end{proposition}

\begin{proof}
We start by applying the law of large numbers result given in (\ref{ML2}) and (\ref{ML3}) to $ \tilde c_N ( \alpha ) := \bar c_N ( \alpha) / \sqrt{ \bar m!} .$
Then $\mathbb{T}( X) = \sigma_N^2 ( \tilde c_N, \varphi ( X) ) , $ if we put $ \varphi( X) = ( \varphi (X^1), \ldots, \varphi ( X^N) ),$ 
and (\ref{ML3})  gives 
\begin{equation*}
\E\left\vert \mathbb{T} (X)-(\E_{t_{k}}^{\Delta }\varphi
(X^{1}))^{\overline{m}}\right\vert ^{2}\leq C_{\bar m} \left\Vert \varphi \right\Vert _{\infty }^{2\overline{m}}\widehat{%
\delta }_{N}^{2}(\tilde {c}_{N})\leq  C_{\bar m} \left\Vert \varphi \right\Vert _{\infty }^{2\overline{m}}N^{-1}.
\end{equation*}

Using now (\ref{pf1}), 
\begin{multline*}
\E\left\vert \E_{t_{k}}(h_{t_k,t_{k+1} }(X)-\widehat{h}_{t_{k},t_{k+1}}(X))\right%
\vert _{N}^{2} \\
=\frac{1}{N}\sum_{i=1}^{N}\E\left\vert
\E_{t_{k}}(H_{t_{k},t_{k+1}}^{i}(X^{i},\mathbb{T}(X))-H_{t_{k},t_{k+1}}^{i}(X^{i},(\E_{t_{k}}^{\Delta }\varphi
(X^{i}))^{\overline{m}})\right\vert ^{2} \\
\leq \frac{C_{stb}^{2}(t_{k+1}-t_{k})^{2}}{N}\sum_{i=1}^{N}\E\left\vert 
\mathbb{T}(X)-(\E_{t_{k}}^{\Delta }\varphi (X^{i}))^{%
\overline{m}}\right\vert ^{2} \\
=C_{stb}^{2}(t_{k+1}-t_{k})^{2}\E\left\vert \mathbb{T} (X)-(\E_{t_{k}}^{\Delta }\varphi (X^{1}))^{\overline{m}%
}\right\vert ^{2}\\
\leq \frac{C_{\bar m} C_{stb}^{2}(t_{k+1}-t_{k})^{2}\left\Vert \varphi \right\Vert _{\infty
}^{2\overline{m}}}{N}.
\end{multline*}%
The argument for $\E\left\vert \widetilde{\E}_{t_{k}}(h_{t_{k},t_{k+1}}(X)-%
\widehat{h}_{t_{k},t_{k+1}}(X))\right\vert ^{2}$ is analogous, based on (\ref%
{pf1'}), and \textbf{B} follows once more from the law of large numbers. 
\end{proof}

We wish to extend the above result to the chains and to iterate them. However, even if 
 $X^{1},\ldots ,X^{N}$ are i.i.d., the coordinates $\widehat{X}_{t_{k}}^{1},...,\widehat{X}%
_{t_{k}}^{N}$ are not independent any more, since they are constructed using all the same 
noise increments $\Delta _{0},...,\Delta
_{k-1}.$ But
freezing $\Delta _{0},...,\Delta _{k-1},$ they become independent. In other 
words, they are independent conditionally on $\Delta
_{0},...,\Delta _{k-1},$ that is, conditionally on $ \mathcal F_{t_k}^\Delta.$ 
Applying once more the conditional law of large numbers (\ref{ML3}) to $ \tilde c_N ( \alpha ) := \bar c_N ( \alpha) / \sqrt{ \bar m!} ,$
we then obtain 
\begin{equation}
\E\left\vert \mathbb{T}(\widehat{X}_{t_{k}})-(\E_{t_{k}}^{%
\Delta }\varphi (\widehat{X}_{t_{k}}^{1}))^{\overline{m}}\right\vert
^{2}\leq C_{\bar m} \left\Vert \varphi \right\Vert
_{\infty }^{2\overline{m}}N^{-1}.  \label{LNN1}
\end{equation}%
and 
\begin{equation}
\E\left[ \left| \sigma _{N}^{2}(\widehat{X}_{t_{k}})-(\E_{t_{k}}^{\Delta }f^2(%
\widehat{X}_{t_{k}}^{1}))^{m} \right|^{2}\right] \leq C_{m} \overline{f}^{4m} N^{ - 1 } .  \label{LNN2}
\end{equation}

Now, in order to compare $\widehat{X}$ with $\overline{X},$ we write%
\begin{equation*}
\widehat{X}_{t_{k+1}}^{i}=\widehat{X}_{t_{k}}^{i}+h_{t_{k},t_{k+1}}^{i}(%
\widehat{X}_{t_{k}})+\bar{S}_{t_{k},t_{k+1}}^{N}(c_N, \widehat{X}%
_{t_{k}})+I_{k}^{i}+J_{k}
\end{equation*}%
where%
\begin{equation}
I_{k}^{i}=\widehat{h}_{t_{k},t_{k+1}}^{i}(\widehat{X}%
_{t_{k}})-h_{t_{k},t_{k+1}}^{i}(\widehat{X}_{t_{k}}),\quad J_{k}=%
\widehat{S}_{t_{k},t_{k+1}}^{N}(c_N, \widehat{X}_{t_{k}})-\bar{S}%
_{t_{k},t_{k+1}}^{N}(c_N, \widehat{X}_{t_{k}} ) .   \label{LNN3}
\end{equation}%
\bigskip

\begin{proposition}
Under the assumptions of Theorem \ref{theo:main},  it holds
that 
\begin{equation}
\E\left\vert \sum_{j=0}^{k}I_{j}\right\vert _{N}^{2}\leq \frac{%
C_{\bar  m}  (1+t_{k+1})t_{k+1}\left\Vert \varphi \right\Vert _{\infty }^{2\overline{m}}}{N}
\label{LN3}
\end{equation}%
and 
\begin{equation}
\E\left\vert \sum_{j=0}^{k}J_{j}\right\vert ^{2}\leq \frac{   C_{ m}  t_{k+1} \underline f^{ - 2 m} \overline f^{4m} }{N}.
\label{LN4}
\end{equation}
\end{proposition}

\begin{proof}
We write%
\begin{equation*}
\E\left\vert \sum_{j=0}^{k}I_{j}\right\vert _{N}^{2}\leq 2\E\left\vert
\sum_{j=0}^{k}\E_{t_{j}}(I_{j})\right\vert _{N}^{2}+2\E\left\vert
\sum_{j=0}^{k}\widetilde{\E}_{t_j}(I_{j})\right\vert _{N}^{2}.
\end{equation*}%
Denote $\gamma _{j}=t_{j+1}-t_{j}.$ Then, 
\begin{equation*}
\E\left\vert \sum_{j=0}^{k}\E_{t_{j}}(I_{j})\right\vert _{N}^{2}=\E\left\vert
\sum_{j=0}^{k}\gamma _{j}\frac{\E_{t_{j}}(I_{j})}{\gamma _{j}}\right\vert
_{N}^{2}\leq t_{k+1 }\sum_{j=0}^{k}\gamma _{j}\E\left\vert \frac{%
\E_{t_{j}}(I_{j})}{\gamma _{j}}\right\vert _{N}^{2} .
\end{equation*}%
Using the Lipschitz property (\ref{pf1}) together with (\ref{LNN1}), we obtain%
\begin{equation*}
\E\left\vert \E_{t_{j}}(I_{j})\right\vert _{N}^{2}\leq C \gamma
_{j}^{2}\E\left\vert \mathbb{T}(\widehat{X}%
_{t_{j}})-(\E_{t_{j}}^{\Delta } \varphi (\widehat{X}^1_{t_{j}}))^{\overline{m}%
}\right\vert ^{2}\leq C^2_{stb}C_{\bar  m}  \left\Vert
\varphi \right\Vert _{\infty }^{2\overline{m}}\gamma _{j}^{2}N^{-1}.
\end{equation*}%
For the last inequality we have used the conditional law of large numbers, based on the fact 
that the components of $\widehat{X}_{t_{j}}$ are independent conditionally
to $\mathcal{F}_{t_{j}}^{\Delta }.$ Putting things together we see that 
\begin{equation*}
\E\left\vert \sum_{j=0}^{k}\E_{t_{j}}(I_{j})\right\vert _{N}^{2}\leq C_{\bar  m}  \left\Vert \varphi \right\Vert _{\infty
}^{2\overline{m}}t_{k+1}^{2}N^{-1}.
\end{equation*}

Using now orthogonality first and then (\ref{pf1'}) together with (\ref{LNN1}), we also obtain
\begin{eqnarray*}
\E\left\vert \sum_{j=0}^{k}\widetilde{\E}_{t_j}(I_{j})\right\vert _{N}^{2}
&=&\sum_{j=0}^{k}\E\left\vert \widetilde{\E}_{t_j}(I_{j})\right\vert
_{N}^{2}\leq C\sum_{j=0}^{k}\gamma _{j}\E\left\vert \mathbb{T}%
(\widehat{X}_{t_{j}})-(\E_{t_{j}}^{\Delta }\varphi(\widehat{X}^1_{t_{j}}))^{\overline{m}}\right\vert ^{2} \\
&\leq &C_{\bar  m}  \left\Vert \varphi
\right\Vert _{\infty }^{2\overline{m}}t_{k+1}N^{-1}.
\end{eqnarray*}

We turn now to the second estimates. The $J_{j}$ are orthogonal. Moreover, by (\ref{CLT4}), we have that  
\begin{equation*}
\E\left\vert J_{j}\right\vert ^{2}= \gamma_j \E\left[  \left| \sigma _{N} (\widehat{X%
}_{t_{j}})- (\E_{t_{j}}^{\Delta }f^2(\widehat{X}_{t_{j}}^{1}))^{m/2} \right|^2 \right]
\leq C_{ m}  \gamma_j  \underline f^{ - 2 m } \overline f ^{4m} N^{ - 1}. 
\end{equation*}%
All in all this gives
\begin{equation*}
\E\left\vert \sum_{j=0}^{k}J_{j}\right\vert ^{2}=\sum_{j=0}^{k}\E\left\vert
J_{j}\right\vert ^{2}\leq C_{m}    t_{k+1} \underline f^{ - 2 m } \overline f ^{4m} N^{-1} . 
\end{equation*}%
\end{proof}

\begin{lemma}
Grant the assumptions of Theorem \ref{theo:main}. Then
\begin{multline}
\E\left\vert \sum_{j=0}^{k}(h_{t_{j},t_{j+1}}(\widehat{X}%
_{t_{j}})-h_{t_{j},t_{j+1}}(\overline{X}_{t_{j}}))\right\vert _{N}^{2}
\label{LN3'} \\
\leq C (1+ L_{\mathbb{T}}^2) (1+t_{k+1})\sum_{j=0}^{k}\E\left\vert \overline{X}_{t_{j}}-%
\widehat{X}_{t_{j}}\right\vert _{N}^{2}(t_{j+1}-t_{j}).
\end{multline}

\end{lemma}

\begin{proof} Let us write $D_{j}=h_{t_{j},t_{j+1}}(\widehat{X}%
_{t_{j}})-h_{t_{j},t_{j+1}}(\overline{X}_{t_{j}}) .$ Then Assumption \ref{ass:H} and (\ref{McKV4}) imply that 
\begin{equation*}
\E\left\vert \E_{t_{j}}(D_{j})\right\vert _{N}^{2}\leq C (1+ L_{\mathbb{T}}^2) \gamma _{j}^{2} \E \left\vert \widehat{X}_{t_{j}}-\overline{X}%
_{t_{j}}\right\vert _{N}^{2}.
\end{equation*}
So
\begin{eqnarray*}
\E\left\vert \sum_{j=0}^{k}\E_{t_{j}}(D_{j})\right\vert _{N}^{2}
&=&\E\left\vert \sum_{j=0}^{k}\gamma _{j}\frac{\E_{t_{j}}(D_{j})}{\gamma _{j}}%
\right\vert _{N}^{2}\leq t_{k+1 }\sum_{j=0}^{k}\gamma _{j}\E\left\vert 
\frac{\E_{t_{j}}(D_{j})}{\gamma _{j}}\right\vert _{N}^{2} \\
&\leq &C t_{k+1} (1+ L_{\mathbb{T}}^2) \sum_{j=0}^{k}\gamma _{j}\E \left\vert \widehat{X}_{t_{j}}-\overline{X}%
_{t_{j}}\right\vert _{N}^{2}.
\end{eqnarray*}%
Moreover,  using orthogonality,
\begin{equation*}
\E\left\vert \sum_{j=0}^{k}\widetilde{\E}_{t_{j}}(D_{j})\right\vert
_{N}^{2}=\sum_{j=0}^{k}\E\left\vert \widetilde{\E}_{t_{j}}(D_{j})\right\vert
_{N}^{2}\leq C (1+ L_{\mathbb{T}}^2)\sum_{j=0}^{k}\gamma
_{j} \E \left\vert \widehat{X}_{t_{j}}-\overline{X}_{t_{j}}\right\vert _{N}^{2},
\end{equation*}%
which allows us to conclude. 
\end{proof}

We are now able to summarize our findings in the following

\begin{theorem}
\label{LLNW2}\textbf{A} Suppose that  (\ref{McKV6}) holds and grant Assumptions \ref{ass:H} and  %
\ref{ass:perturbation} for all $ (s, t) = (t_i, t_{i+1}), $ such that $ t_{i+1} \le T.$ 
Suppose furthermore that  for some $ q \geq 3,$ $ \gamma_T( q) < \infty $ (recall (\ref{eq:gammaTq})).

Let $X_{0} , \widehat X_0 $ and $ \overline X_0 \in L^2_0 ( \R^N) $ be composed of  independent and identically distributed random variables each.

Then for any $k \geq  0 $
such that $ t_{k+1} \le T , $  
\begin{equation}
\E\left\vert \overline{X}_{t_{k+1}}-\widehat{X}_{t_{k+1}} \right\vert
_{N}^{2}\leq 5e^{Q_{1}(t_{k+1}+1)^{2}}\left\{ \E\left\vert \overline{X}%
_{{0}}^{N}-\widehat{X}_{{0}}^{N}\right\vert _{N}^{2}+Q_{2}\frac{1}{N}%
t_{k+1} \right\} ,  \label{clt5}
\end{equation}%
where $Q_{1}$ and $Q_{2}$ are constants depending on $T,C_{stb},L_{\varphi },$ $\overline{L},$ $m, \bar m, \| \varphi\|_\infty, \underline f, \bar f$  only.%

\textbf{B }Moreover, for every $p<p(q)$ (recall (\ref{pq})), there exists a constant $Q $ 
depending only on $m, \bar m, p, q, \gamma_T(q),$ on $ L_\varphi , L_{f^2 } , C_{stb},$ on $\|\varphi\|_\infty, \underline f , \overline f, $ and on the constant $C$ of 
(\ref{McKV6}), 
such that 
\begin{equation}
\mathcal{W}_{2}^{2}(X_{t_{k}},\widehat{X}_{t_{k}})\leq Q e^{QT^2 }\left( 
\mathcal{W}_{2}^{2}(X_{0},\widehat{X}_{0})+t_{k}(\frac{1}{N^{p/2}}%
+\sum_{r=1}^{m-1}\left\Vert c_{N}\otimes _{r}c_{N}\right\Vert
_{m-r,2}^{2})\right) .  \label{CLT6}
\end{equation}%
If Assumption \ref{ass:2.10first} holds for all $ q \geq 3, $ then the above estimate holds for every $p<\frac{1}{3}.$
\end{theorem}

\begin{proof}
\textbf{A } We use the optimal coupling of $ (\bar X_0, \widehat X_0 ),$ for the Wasserstein-$2-$distance $\mathcal{W}_{2}.$ 
Since $\bar X_0$ and $ \widehat X_0$ have i.i.d. components, this optimal coupling still produces i.i.d. components (of optimally coupled pairs). We construct this coupling such that it is independent of all innovations. We write 
\begin{eqnarray*}
\widehat{X}_{t_{k+1}}^{i}-\overline{X}_{t_{k+1}}^{i} &=&\widehat{X}%
_{{0}}^{i}-\overline{X}_{{0}}^{i}+\sum_{j=0}^{k}(h_{t_{j},t_{j+1}}^{i}(%
\widehat{X}_{t_{j}})-h_{t_{j},t_{j+1}}^{i}(\overline{X}_{t_{j}})) \\
&&+\sum_{j=0}^{k}(\bar{S}_{t_{j},t_{j+1}}^{N}(c_N, \widehat{X}_{t_{j}}^{N})-\bar{S}%
_{t_{j},t_{j+1}}^{N}(c_N, \overline{X}_{t_{j}}^{N}))+\sum_{j=0}^{k}I_{j}^{i}+%
\sum_{j=0}^{k}J_{j},
\end{eqnarray*}%
where we use the same notation as in (\ref{LNN3}). We have proved in (\ref%
{LN3'}) that 
$$
\E\left\vert \sum_{j=0}^{k}(h_{t_{j},t_{j+1}}(\widehat{X}%
_{t_{j}})-h_{t_{j},t_{j+1}}(\overline{X}_{t_{j}}))\right\vert _{N}^{2} 
\leq C (1+ L_{\mathbb{T}}^2) (1+t_{k+1}) \sum_{j=0}^{k}\E\left\vert \overline{X}_{t_{j}}-%
\widehat{X}_{t_{j}}\right\vert _{N}^{2}(t_{j+1}-t_{j}).
$$

Moreover, the $\bar{S}_{t_{j},t_{j+1}}^{N}(c_N, \widehat{X}_{t_{j}})-\bar{S}%
_{t_{j},t_{j+1}}^{N}(c_N, \overline{X}_{t_{j}}),j=0,...,k,$ are orthogonal,
therefore, recalling also (\ref{McKV9'}), 
\begin{multline*}
\E\left\vert \sum_{j=0}^{k}\bar{S}_{t_{j},t_{j+1}}^{N}(c_N, \widehat{X}%
_{t_{j}}^{N})-\bar{S}_{t_{j},t_{j+1}}^{N}(c_N, \overline{X}_{t_{j}}^{N})\right%
\vert ^{2} \\
=\sum_{j=0}^{k}\E\left\vert \bar{S}_{t_{j},t_{j+1}}^{N}(c_N, \widehat{X}%
_{t_{j}}^{N})-\bar{S}_{t_{j},t_{j+1}}^{N}(c_N, \overline{X}_{t_{j}}^{N})\right%
\vert ^{2}\leq \overline{L}\sum_{j=0}^{k}\E\left\vert \widehat{X}_{t_{j}}^{N}-%
\overline{X}_{t_{j}}^{N}\right\vert _{N}^{2}(t_{j+1}-t_{j}).
\end{multline*}%
Recall (\ref{LN3}) and (\ref{LN4}).\ Gathering these results, we obtain that 
\begin{multline*}
\E\left\vert \overline{X}_{t_{k+1}}-\widehat{X}_{t_{k+1}}\right\vert
_{N}^{2}\\
\leq 5 \E\left\vert \overline{X}_{t_{0}}-\widehat{X}%
_{t_{0}}\right\vert _{N}^{2}+C_{1}\sum_{j=0}^{k}\E\left\vert \overline{X}%
_{t_{j}}-\widehat{X}_{t_{j}}\right\vert
_{N}^{2}(t_{j+1}-t_{j})+C_{2}(t_{k+1}-t_{0})N^{-1},
\end{multline*}%
for some constants $C_{1},C_{2}$ depending on $T,$ and on all model parameters. Using the discrete version of the Gronwall lemma allows
to conclude.

In order to obtain \textbf{B}, we take $ \bar X_0 := X_0 $ and we use once more the optimal coupling of $ (X_0, \widehat X_0 ),$ for the Wasserstein-$2-$distance $\mathcal{W}_{2},$ such that it is independent of all innovations. 
Then we use that $ \mathcal{W}_{2}^{2} (X_{t_{k}},\widehat{X}_{t_{k}}) \le 2  \mathcal{W}_{2}^{2}( X_{t_{k}},\overline {X}_{t_{k}}) + 2  \mathcal{W}_{2}^{2} (\overline X_{t_{k}},\widehat{X}_{t_{k}} )$
and apply the coupling estimate obtained in (\ref{CLT5}).
\end{proof}

\begin{proof}[Proof of Theorem \ref{theo:main}]
Theorem \ref{theo:main} follows now directly from Theorem \ref{LLNW2}, by taking $ X_0 = \widehat X_0 .$ 

\end{proof}

\section{An example: the quadratic case}\label{sec:5}
In this section we present an example which corresponds to the so-called
``quadratic Central Limit Theorem''. Important previous papers on this subject are
\cite{R75}, \cite{R79}, \cite{dJ89}, \cite{dJ90}, \cite{G91}, \cite{RR96} and \cite{DR09}.

We restrict ourselves to the case $m=2.$ Similar examples can be given for multi-linear homogeneous sums of higher order. However, the computations then become more involved.

\subsection{The coefficients and the error bound}
Throughout this section we suppose that $N\geq 2.$ We fix $\theta\geq \frac{1}{2}$ and we consider the coefficients
\begin{equation}\label{eq:cthetaN}
c_{\theta,N}(i,j)=\kappa_{N}(\theta )\,\alpha_N(\theta)\,\left\vert j-i\right\vert ^{-\theta }\mathbf{1}_{\{i\neq j\}},\quad 1\le i,j\le N.
\end{equation}
Here
\begin{equation*}
\alpha_{N}(\theta)=\left\{
\begin{array}{ll}
\frac{1}{2\sqrt{J(2\theta )N}}, & \theta>\frac{1}{2}, \\[1mm]
\frac{1}{2\sqrt{N\log N}}, & \theta=\frac{1}{2},
\end{array}
\right.
\qquad \mbox{where}\quad J(s):=\sum_{k=1}^{\infty }k^{-s},\quad s>1.
\end{equation*}
The constant $\kappa_N(\theta)>0$ is chosen such that the normalization condition (\ref{eq:rencN}) holds. Since $m=2,$ this condition reads
$$\left\Vert c_{\theta,N}\right\Vert _{2}^{2}=\frac{1}{2}.$$

With this choice of coefficients,  we have the following result. 
\begin{theorem}\label{theo:example}
Let $\theta\geq\frac{1}{2},$ and let $c_N=c_{\theta,N}$ be given by (\ref{eq:cthetaN}). Moreover we suppose that the assumptions of Theorem \ref{theo:main}
hold for all $ q > 3.$ 

Then for every $p<\frac{1}{3}$ there exists a constant $Q,$ depending on $\theta,$ $p,$ and on the model parameters, such that for every $t_k\leq T,$
\begin{equation}
\mathcal{W}_{2}^{2}(X_{t_{k}},\widehat{X}_{t_{k}})\leq Qe^{QT^{2}}\,r_{N,p}(\theta),\quad
r_{N,p}(\theta):=\left\{
\begin{array}{ll}
\frac{1}{(\log N)^{2}}, & \theta=\frac{1}{2}, \\
N^{-p/2}+N^{-(4\theta -2)}, & \frac{1}{2}<\theta<\frac{3}{4}, \\
N^{-p/2}, & \theta\geq\frac{3}{4}.
\end{array}
\right.  \label{eq:rateexample}
\end{equation}
\end{theorem}

The remainder of this section is devoted to the proof of Theorem \ref{theo:example} 

\subsection{Proof of Theorem \ref{theo:example}}
We start with some elementary facts.  We write $C_\theta$ for a constant which depends only on $\theta.$ Its value may change from line to line. For $M\geq 1$ and $s\geq 0$ we put
$$H_{M}^{(s)}=\sum_{k=1}^{M}k^{-s}.$$
Comparing sums with integrals, we have the classical asymptotics
\begin{eqnarray}
H_{M}^{(s)} &=&J(s)-\frac{M^{1-s}}{s-1}+O(M^{-s}),\quad s>1,  \label{eq:H1} \\
H_{M}^{(s)} &=&\frac{M^{1-s}}{1-s}+O(1),\quad 0<s<1,  \label{eq:H2} \\
H_{M}^{(1)} &=&\log M+\gamma +\frac{1}{2M}+O(M^{-2}),  \label{eq:H3}
\end{eqnarray}%
where $\gamma $ is Euler's constant. We shall also use the following elementary bounds, valid for every $M\geq 1:$
\begin{equation}
\begin{gathered}
H_{M}^{(s)}\leq \frac{M^{1-s}}{1-s}\;\;(0<s<1),\qquad H_{M}^{(1)}\leq 1+\log M,\\
H_{M}^{(s)}\leq J(s),\qquad \sum_{k>M}k^{-s}\leq \frac{M^{1-s}}{s-1}\;\;(s>1).
\end{gathered}  \label{eq:Hbounds}
\end{equation}

\subsubsection{Normalization}
We first check that the constant $\kappa_N(\theta)$ is bounded and that it converges to $1.$

\begin{lemma}[Normalization]\label{lem:normalization}
We have
\begin{equation*}
\left\Vert c_{\theta,N}\right\Vert _{2}^{2}=\kappa_N^2(\theta) \left( \frac{1}{2}+O(\eta_N(\theta))\right) ,
\quad \mbox{where}\quad
\eta_N(\theta)=\left\{
\begin{array}{ll}
\frac{1}{\log N}, & \theta=\frac{1}{2}, \\
N^{-(2\theta-1)}, & \frac{1}{2}<\theta<1, \\
\frac{\log N}{N}, & \theta=1, \\
\frac{1}{N}, & \theta>1.
\end{array}
\right.
\end{equation*}
In particular, $\kappa_N(\theta)\to 1$ as $N\to\infty,$ and $\sup_{N\geq 2}\kappa_N(\theta)<\infty.$
\end{lemma}

\begin{proof}
We have that 
$$\left\Vert c_{\theta,N}\right\Vert _{2}^{2}=2\kappa_N^2(\theta)\alpha_N^2(\theta)\sum_{1\le i<j\le N}|j-i|^{-2\theta}.$$
For every $1\leq k\leq N-1,$ there are exactly $N-k$ couples $i<j$ such that $j-i=k.$ Hence
\begin{equation}
\sum_{1\leq i<j\leq N}\left\vert j-i\right\vert ^{-2\theta }=\sum_{k=1}^{N-1}%
\frac{N-k}{k^{2\theta }}=NH_{N-1}^{(2\theta )}-H_{N-1}^{(2\theta -1)}.
\label{eq:sumpairs}
\end{equation}

\textbf{Case $\theta>\frac{1}{2}$.} In this case $2\alpha_N^2(\theta)=\frac{1}{2J(2\theta)N}.$ So it is sufficient to show that
\begin{equation}\label{eq:normgoal}
\frac{1}{N}\sum_{1\leq i<j\leq N}\left\vert j-i\right\vert ^{-2\theta }=J(2\theta)+O(\eta_N(\theta)).
\end{equation}
We use (\ref{eq:H1}) with $s=2\theta$ and $M=N-1.$ This gives
$$NH_{N-1}^{(2\theta )}=NJ(2\theta)-\frac{N^{2-2\theta}}{2\theta-1}+O(N^{1-2\theta}).$$
We now distinguish three sub-cases.

\emph{Sub-case $\frac{1}{2}<\theta<1$.} By (\ref{eq:H2}) with $s=2\theta-1,$
$$H_{N-1}^{(2\theta -1)}=\frac{N^{2-2\theta }}{2-2\theta }+O(1).$$
Therefore
\begin{equation*}
\frac{1}{N}\sum_{1\leq i<j\leq N}\left\vert j-i\right\vert ^{-2\theta }=J(2\theta )-\frac{N^{1-2\theta}}{(2\theta-1)(2-2\theta)}+O\Big(\frac{1}{N}\Big)=J(2\theta)+O(N^{-(2\theta-1)}).
\end{equation*}

\emph{Sub-case $\theta=1$.} Here $NH_{N-1}^{(2)}=NJ(2)-1+O(N^{-1}).$ Moreover, by (\ref{eq:H3}), $H_{N-1}^{(1)}=\log N+\gamma +O(N^{-1}).$ Hence
\begin{equation*}
\frac{1}{N}\sum_{1\leq i<j\leq N}\left\vert j-i\right\vert ^{-2}=J(2)-\frac{\log N+\gamma +1}{N}+O(N^{-2})=J(2)+O\Big(\frac{\log N}{N}\Big).
\end{equation*}

\emph{Sub-case $\theta>1$.} Since $2-2\theta<0,$ we have $NH_{N-1}^{(2\theta )}=NJ(2\theta)+O(1).$ Moreover, $0\leq H_{N-1}^{(2\theta -1)}\leq J(2\theta-1)$ by (\ref{eq:Hbounds}). Hence (\ref{eq:normgoal}) holds with $\eta_N(\theta)=\frac{1}{N}.$

\textbf{Case $\theta=\frac{1}{2}$.} By (\ref{eq:sumpairs}),
$$\sum_{1\leq i<j\leq N}\left\vert j-i\right\vert ^{-1}=NH_{N-1}^{(1)}-(N-1).$$
By (\ref{eq:H3}), $H_{N-1}^{(1)}=H_{N}^{(1)}-\frac{1}{N}=\log N+\gamma -\frac{1}{2N}+O(N^{-2}).$ We obtain
\begin{equation*}
\sum_{1\leq i<j\leq N}\left\vert j-i\right\vert ^{-1}=N\log N+(\gamma -1)N+O(1).
\end{equation*}
Since $2\alpha_N^2(\frac{1}{2})=\frac{1}{2N\log N},$ this gives
$$\left\Vert c_{\frac{1}{2},N}\right\Vert _{2}^{2}=\frac{\kappa_N^2(\frac12)}{2}\left(1+\frac{\gamma-1}{\log N}+O\Big(\frac{1}{N\log N}\Big)\right).$$

We are now able to conclude our proof. Since $\left\Vert c_{\theta,N}\right\Vert _{2}^{2}=\frac{1}{2},$ we have $\kappa_N^2(\theta)=(1+O(\eta_N(\theta)))^{-1}.$ Hence $\kappa_N(\theta)\to 1.$ Moreover, $\kappa_N(\theta)<\infty$ for each $N.$ Therefore $\sup_{N\geq 2}\kappa_N(\theta)<\infty.$
\end{proof}

\subsubsection{Influence and contraction factors}
We now compute the two quantities which enter the error bound of Theorem \ref{theo:main}: the influence factor and the contraction factor. The contraction factor involves discrete convolutions of the kernel $k\mapsto |k|^{-\theta}.$ We start with an elementary estimate on these convolutions.

\begin{lemma}\label{lem:convolution}
For $1\leq i<j\leq N,$ put $r=j-i$ and
\begin{equation*}
A_{i,j}^{(\theta)}=\sum_{k\in\{1,\ldots,N\}\setminus\{i,j\}}\left\vert i-k\right\vert ^{-\theta }\left\vert j-k\right\vert ^{-\theta }.
\end{equation*}
Then $A_{i,j}^{(\theta)}\leq C_\theta\, a_{N,\theta}(r),$ where
\begin{equation*}
a_{N,\theta}(r)=\left\{
\begin{array}{ll}
1+\log \frac{N}{r}, & \theta=\frac{1}{2}, \\
r^{1-2\theta }, & \frac{1}{2}<\theta<1, \\
\frac{1+\log r}{r}, & \theta=1, \\
r^{-\theta }, & \theta>1.
\end{array}
\right.
\end{equation*}
\end{lemma}

\begin{proof}
We put $u=k-i $ such that $ u \in \{-(N-1),\ldots ,N-1\}\setminus \{0,r\}.$ We split the sum into three parts, according to $u<0,$ $0<u<r$ and $u>r.$

For $u<0$ we put $w=-u.$ For $u>r$ we put $w=u-r.$ In both cases $1\leq w\leq N-1,$ and the corresponding term equals $w^{-\theta }(w+r)^{-\theta }.$ Therefore
$$A_{i,j}^{(\theta)}\leq 2B_{N}(r)+M(r),$$
where
\begin{equation*}
B_{N}(r):=\sum_{w=1}^{N}w^{-\theta }(w+r)^{-\theta },\qquad
M(r):=\sum_{u=1}^{r-1}u^{-\theta }(r-u)^{-\theta },
\end{equation*}
with the convention $M(1)=0.$ We bound the two terms separately.

\textbf{The term $M(r)$.} We use the symmetry of $M(r) $ with respect to $u\mapsto r-u.$ Since $r-u\geq r/2$ for $u\leq r/2,$ we see that 
\begin{equation*}
M(r)\leq 2\sum_{u \geq 1 : u \leq r/2}u^{-\theta }(r-u)^{-\theta }\leq 2\Big(\frac{r}{2}\Big)^{-\theta }\sum_{u \geq 1 : u \leq r/2}u^{-\theta }.
\end{equation*}
By (\ref{eq:Hbounds}), we deduce that
\begin{equation*}
M(r)\leq \left\{
\begin{array}{ll}
\frac{2^{2\theta }}{1-\theta }\,r^{1-2\theta }, & \theta<1, \\[1mm]
\frac{4(1+\log (r/2) )}{r}, & \theta=1, \\[1mm]
2^{1+\theta }J(\theta )\,r^{-\theta }, & \theta>1.
\end{array}
\right.
\end{equation*}
In what follows we will simply upper bound further $ \log ( r/2 ) \le \log r  $ to obtain the bound $M(r) \le  \frac{4(1+\log r )}{r} $ in case $\theta = 1.$ 
In particular, we notice that $M(r)\leq 4$ in case that  $\theta=\frac{1}{2}.$

\textbf{The term $B_N(r)$.} We split the sum according to $w\leq r$ and $w>r.$ For $w\leq r$ we use $(w+r)^{-\theta }\leq r^{-\theta }.$ For $w>r$ we use $(w+r)^{-\theta }\leq w^{-\theta }.$ Together with (\ref{eq:Hbounds}), this gives the following bounds.

If $\theta=\frac{1}{2},$ then
\begin{equation*}
B_{N}(r)\leq r^{-1/2}\sum_{w=1}^{r}w^{-1/2}+\sum_{w=r+1}^{N}\frac{1}{w}\leq 2+\int_{r}^{N}\frac{dx}{x}=2+\log \frac{N}{r}.
\end{equation*}

If $\frac{1}{2}<\theta<1,$ then
\begin{equation*}
B_{N}(r)\leq r^{-\theta }\frac{r^{1-\theta }}{1-\theta }+\sum_{w>r}w^{-2\theta }\leq \Big(\frac{1}{1-\theta }+\frac{1}{2\theta -1}\Big)r^{1-2\theta }.
\end{equation*}

If $\theta=1,$ then
\begin{equation*}
B_{N}(r)\leq \frac{1+\log r}{r}+\sum_{w>r}w^{-2}\leq \frac{2+\log r}{r}.
\end{equation*}

If $\theta>1,$ then, since $1-2\theta\leq -\theta$ and since $ r \geq 1, $
\begin{equation*}
B_{N}(r)\leq J(\theta )r^{-\theta }+\frac{r^{1-2\theta }}{2\theta -1}\leq \Big(J(\theta )+\frac{1}{2\theta -1}\Big)r^{-\theta }.
\end{equation*}

Putting  these estimates together concludes the proof.
\end{proof}

\begin{lemma}\label{lem:factors}
For every $\theta\geq\frac{1}{2}$ and every $N\geq 2,$
\begin{equation}
\delta _{N}(c_{\theta,N})+\widehat{\delta }_{N}(c_{\theta,N})\leq \frac{C_\theta}{\sqrt{N}}  \label{eq:influencetheta}
\end{equation}
and
\begin{equation}
\left\Vert c_{\theta,N}\otimes _{1}c_{\theta,N}\right\Vert _{1,2}^{2}\leq C_\theta\, \rho_N(\theta),\quad \mbox{where}\quad
\rho_N(\theta)=\left\{
\begin{array}{ll}
\frac{1}{(\log N)^{2}}, & \theta=\frac{1}{2}, \\
N^{-(4\theta -2)}, & \frac{1}{2}<\theta<\frac{3}{4}, \\
\frac{1+\log N}{N}, & \theta=\frac{3}{4}, \\
\frac{1}{N}, & \theta>\frac{3}{4}.
\end{array}
\right.  \label{eq:contractiontheta}
\end{equation}
\end{lemma}

\begin{proof}
We write $c=c_{\theta,N},$ $\kappa=\kappa_N(\theta)$ and $\alpha=\alpha_N(\theta)$ for short. 

\textbf{Influence factors.} Since $m=2,$ the multi-indices of $\Lambda_{N,2}^{i}$ are $(i,j)$ and $(j,i),$ with $j\neq i.$ Therefore
\begin{equation*}
\delta _{N}^{2}(c)=2\max_{i}\sum_{j=1, j \neq i }^{N}c^{2}(i,j)
\end{equation*}
and
\begin{equation*}
\widehat{\delta }_{N}^{2}(c)=\sum_{i=1}^{N}\Big(\sum_{j=1, j \neq i }^{N}c^{2}(i,j)\Big)^{2}\leq N\max_{i}\Big(\sum_{j=1, j \neq i }^{N}c^{2}(i,j)\Big)^{2}.
\end{equation*}
For any fixed $i, $ we have by (\ref{eq:Hbounds}),
\begin{equation*}
\sum_{j: j \neq i}|j-i|^{-2\theta}\leq 2H_{N}^{(2\theta)}\leq \left\{
\begin{array}{ll}
2J(2\theta), & \theta>\frac{1}{2}, \\
2(1+\log N), & \theta=\frac{1}{2}.
\end{array}
\right.
\end{equation*}
Recalling the definition of $\alpha_N(\theta),$ we obtain in both cases
\begin{equation}
\max_{i}\sum_{j=1}^{N}c^{2}(i,j)\leq \frac{C_\theta}{N}, \label{eq:rowsum}
\end{equation}
which implies (\ref{eq:influencetheta}).

\textbf{Contraction factor.} By (\ref{inf4}) with $m=2$ and $q=1,$
\begin{equation*}
\left\Vert c\otimes _{1}c\right\Vert _{1,2}^{2}=\sum_{i,j=1}^{N}\Big(\sum_{k=1}^{N}c(i,k)c(j,k)\Big)^{2}=D_N+O_N,
\end{equation*}
where $D_N$ is the sum over the diagonal elements $i=j,$ and $O_N$ the sum over $i\neq j.$

\emph{The diagonal part.} By (\ref{eq:rowsum}),
\begin{equation*}
D_N=\sum_{i=1}^{N}\Big(\sum_{k=1}^{N}c^{2}(i,k)\Big)^{2}\leq N\frac{C_\theta}{N^{2}}=\frac{C_\theta}{N}.
\end{equation*}

\emph{The off-diagonal part.} For $i\neq j,$ we have $\sum_{k}c(i,k)c(j,k)=\kappa^{2}\alpha^{2}A_{i\wedge j,i\vee j}^{(\theta)}.$ Moreover, there are $2(N-r)$ couples $(i,j)$ such that $|i-j|=r.$ Lemma \ref{lem:convolution} thus gives
\begin{equation}
O_N\leq C_\theta\,\alpha^{4}\sum_{r=1}^{N-1}(N-r)a_{N,\theta}^{2}(r)\leq C_\theta\,\alpha^{4}N\sum_{r=1}^{N-1}a_{N,\theta}^{2}(r).  \label{eq:ON}
\end{equation}
We now study the right hand side of (\ref{eq:ON}).

If $\theta=\frac{1}{2},$ then $\alpha^{4}=\frac{1}{16N^{2}(\log N)^{2}}.$ The function $x\mapsto (1+\log \frac{1}{x})^{2}$ is decreasing on $(0,1].$ Hence
\begin{equation*}
\sum_{r=1}^{N-1}\Big(1+\log \frac{N}{r}\Big)^{2}\leq N\int_{0}^{1}\Big(1+\log \frac{1}{x}\Big)^{2}dx=N\int_{0}^{\infty }(1+y)^{2}e^{-y}dy=5N.
\end{equation*}
Consequently, $O_N\leq \frac{C}{(\log N)^{2}}.$

If $\theta>\frac{1}{2},$ then $\alpha^{4}=\frac{1}{16J(2\theta)^{2}N^{2}}.$ In this case (\ref{eq:ON}) becomes
$$O_N\leq \frac{C_\theta}{N}\sum_{r=1}^{N-1}a_{N,\theta}^{2}(r).$$
We distinguish the following cases.
\begin{enumerate}
\item If $\frac{1}{2}<\theta<\frac{3}{4},$ then $2-4\theta>-1.$ So $\sum_{r=1}^{N-1}r^{2-4\theta}\leq C_\theta N^{3-4\theta},$ and $O_N\leq C_\theta N^{-(4\theta-2)}.$
\item If $\theta=\frac{3}{4},$ then $\sum_{r=1}^{N-1}r^{2-4\theta} = \sum_{r=1}^{N-1}r^{-1}\leq 1+\log N,$ and $O_N\leq C\frac{1+\log N}{N}.$
\item If $\frac{3}{4}<\theta<1,$ then $\sum_{r\geq 1}r^{2-4\theta}=J(4\theta-2)<\infty,$ and $O_N\leq \frac{C_\theta}{N}.$
\item If $\theta=1,$ then $\sum_{r\geq 1}\frac{(1+\log r)^{2}}{r^{2}}<\infty,$ and $O_N\leq \frac{C}{N}.$
\item If $\theta>1,$ then $\sum_{r\geq 1}r^{-2\theta}<\infty,$ and $O_N\leq \frac{C_\theta}{N}.$
\end{enumerate}

It follows that in all cases, the diagonal part $D_N$ is of lower order $ C_\theta / N $ and thus absorbed. From this we deduce (\ref{eq:contractiontheta}).
\end{proof}

\subsubsection{The proof of Theorem \ref{theo:example}}
The statement of Theorem \ref{theo:example} is now an immediate consequence of Theorem \ref{theo:main} and of Lemma \ref{lem:factors}.

\begin{remark}
For $\frac{1}{2}<\theta<\frac{3}{4},$ the rate in (\ref{eq:rateexample}) is $N^{-\min(p/2,\,4\theta-2)},$ up to a factor $2.$ Since $p<\frac{1}{3},$ the contraction term dominates only for $\theta<\frac{1}{2}+\frac{p}{8}<\frac{13}{24}.$ For larger $\theta,$ the rate is limited by the term $N^{-p/2}.$ This term comes from the passage from $d_3$ to $W_2.$
\end{remark}

\begin{remark}\label{rem:thetasmall}
Let now $0\leq \theta<\frac{1}{2}.$ In this case the natural normalization is $\alpha_N(\theta)\sim C_\theta N^{\theta-1},$ since
$$\sum_{i\neq j}|j-i|^{-2\theta}\sim N^{2-2\theta}\int_{[0,1]^2}|s-t|^{-2\theta}ds\,dt.$$

The influence factor still vanishes. Indeed, $\delta_N^2(c_{\theta,N})\leq C_\theta N^{2\theta-2}H_N^{(2\theta)}\leq \frac{C_\theta}{N}.$ Hence the invariance principle applies.

However, the Gaussian limit fails. In the representation (\ref{CLT9}), the kernel $\phi_N$ converges in $L^2([0,1]^2)$ to $\phi_\theta(s,t)=C_\theta|s-t|^{-\theta}.$ Hence $S_N(c_{\theta,N},G)$ converges to the double Wiener integral $I_2(\phi_\theta),$ see \cite{BC}, Theorem 4.4. Since $I_2(\phi_\theta)$ is a non-zero element of the second Wiener chaos, it is not Gaussian. Equivalently, by the fourth moment theorem \cite{NuPe}, the contraction $\Vert c_{\theta,N}\otimes_1 c_{\theta,N}\Vert_{1,2}$ does not converge to $0.$ For $\theta=0,$ one obtains a multiple of $W_1^2-1,$ that is, a centred chi-square variable.

Our method does not apply in this case. One expects that this gives rise to a non-Brownian common noise in the discrete time limit scheme $\widehat X.$ 
\end{remark}

{\noindent\textbf{Acknowledgements}:
Eva L\"ocherbach has received funding from the European Union's Horizon Europe research and innovation programme under the Marie Sklodowska-Curie Actions Staff Exchanges (Grant Agreement No. 101183168, Call: HORIZON-MSCA-2023-SE-01). 
}

{\noindent\textbf{Disclaimer}: Funded by the European Union. Views and opinions expressed are however those of the author(s) only and do not necessarily reflect those of the European Union or the European Education and Culture Executive Agency (EACEA). Neither the European Union nor EACEA can be held responsible for them.}

\end{document}